\documentclass[preprint,1p]{elsarticle}

\usepackage{amsthm,amsfonts,amsmath,amssymb,empheq}
\usepackage{color}
\definecolor{Blue}{RGB}{0,0,255}
\usepackage[unicode, colorlinks, urlcolor=Blue, linkcolor=Blue, citecolor=Blue]{hyperref}
\usepackage{subcaption}

\usepackage[ruled,vlined]{algorithm2e}
\usepackage{textcomp}
\usepackage[normalem]{ulem}
\usepackage{listings}
\usepackage{multirow}
\newcommand{\sign}{\operatorname{sign}}

\journal{Applied Mathematics and Computation}
 
\usepackage{geometry}
\newtheorem{theorem}{Theorem}[section]

\newtheorem{remark}[theorem]{Remark}
\theoremstyle{definition}

\begin{document}

\begin{frontmatter}

\title{An Algorithm for Solving the Burgers Equation Based on Rational Approximations and the Collocation Method}

\author[2,3]{I.~A.~Bugoets}
\ead{i.bugoets@g.nsu.ru}

\author[1,2,3,4]{B.~V. ~Semisalov\corref{mycorrespondingauthor}}
\cortext[mycorrespondingauthor]{Principal corresponding author and corresponding author}
\ead{vibis87@gmail.com}

\address[1]{Faculty of Computational Mathematics and Cybernetics, Shenzhen MSU-BIT University, Shenzhen, China}
\address[2]{Sobolev Institute of Mathematics, Novosibirsk, Russia}
\address[3]{Novosibirsk State University, Novosibirsk, Russia}
\address[4]{Federal Research Center for Information and Computational Technologies, Novosibirsk, Russia}

\begin{abstract}
 We propose a numerical method for solving initial-boundary value problems for the Burgers equation that takes into account the presence of singularities of the analytic continuation of the solution into the complex plane. The algorithm uses rational approximations, a collocation method and a new scheme for time discretization combined with the relaxation method. An adaptive strategy is developed to  account for the location of the singular point of the solution.
 Solutions to the well-known benchmark problem with a sine-shaped initial condition are obtained for the values of the viscosity as low as $\nu = 10^{-8}$. Verification of the numerical results in such an extreme parameter regime makes it necessary to modify a classical formula for the exact solution so as to avoid floating-point overflow during its computation. The developed method is compared with an algorithm employing the Runge--Kutta scheme, and the conditions under which each approach yields superior results are outlined.
An important feature of the proposed method is its capability to track the singular point of the solution. To verify this, we consider a problem whose solution has a simple analytical form and demonstrate that the exact and approximate coordinates of the singular point of the solution agree with high accuracy, even when the singularity comes extremely close to the domain of the problem. Taking the location of the singular point into account yields significant gains in both accuracy and run time relative to other methods.
\end{abstract}


\end{frontmatter}

\section{Overview of the problem and approaches to solving it}
\label{S1}
Numerical analysis of the non-linear differential equations that describe complex processes of hydro- and gas dynamics presents an essential interest both for the fundamental science and for many applications. One of the most universal benchmarks in this context is solving the Burgers Equation (BE) \cite{Burgers}\footnote{For the first time this equation appeared in Bateman's paper \cite{Bat} several decades before the famous work of Burgers which had initiated its systematic studies.} which is commonly used to test the accuracy, stability and convergence of numerical methods:
\begin{equation}
\label{BE}
u_t + u u_x = \nu u_{xx}.
\end{equation}
Here, $u(t,x)$ is the velocity of flow of liquid, $t$ is the time (in what follows we assume that $t \in [T_0,T]$), $x\in[-1,1]$ is the spatial coordinate, $\nu$ is the kinematic viscosity. 

A specific interest to the BE from the side of specialists in the field of hydrodynamics is related with its considerable simplicity and the presence of two important quantities, $uu_x$ and $\nu u_{xx}$, which characterize the essence of non-linear dynamics of the viscous liquid and present the convective and the viscous terms respectively. Thanks to the presence of these terms, BE describes the basic mechanisms of the emergence and evolution of wave fronts in the incompressible viscous liquid and the pathways to turbulence. Moreover, analysis of the BE enables one to describe the change in the structure of flow while the viscosity is reduced: the emergence of the steep gradients of velocity which causes the breaking (loss of stability) of the wave and other specific phenomena. 

An important circumstance is that an exact solution to the BE is known. The work \cite{ExSols} describes more that 30 different solutions of the BE of the relatively simple forms and proposes the way to classify them. In general case, the exact solution of \eqref{BE} can be obtained using the Cole--Hopf transformation \cite{Hopf, Cole}, which reduces \eqref{BE} to the heat conduction equation. It should be noted, however, that for small values of the viscosity $\nu$, the direct usage of the known analytical formulas to compute exact solutions of the BE in floating point arithmetic may present considerable difficulties, which have not been discussed in the literature yet. These difficulties and the ways to overcome them are described below in section~\ref{S2}.

Many different numerical methods were proposed for solving the BE. The first were the finite-difference schemes \cite{Crank, Lax} which subsequently played an essential role in formation and development of the approaches for simulating a number of different physical processes. Later the schemes of high order of approximation with enhanced properties of the stability and monotonicity were proposed \cite{Hassanien, Zhanlav}. Further in this direction we should mention the application of different kinds of splines for approximating the solution with respect to spatial variable in combination with Crank--Nicolson and higher order schemes for time-stepping \cite{Ali,Sayed,Kutulay}. The finite element methods \cite{Dogan, Ucar} and the spectral methods \cite{Maerschalck, Basdevan} were applied as well. For enhancing the accuracy and stability of methods the following approaches were used: integral representation of the solution \cite{Somrath}; reduction of the equation to a system of first-order differential equations and application of Pade approximations for computing the matrix exponential \cite{Kemal}; approximation of the solution by different types of orthogonal polynomials \cite{Khater,Adel,Arar}. We should also mention the considerably high flexibility and appropriate accuracy of the method based on multilevel wavelets \cite{OVasilyev}. 

The mentioned approaches enabled one to obtain an increased accuracy while solving BE with small viscosity down to values of the order $\nu\sim10^{-4}$ with representative order of the relative error equal to $10^{-2}$--$10^{-3}$ in supremum norm, which is good enough for engineering simulations. However, further reduction of the viscosity and of the error plays an important role for fundamental research, in particular, for the analysis of hydrodynamical equations with very small values of $\nu$ including the case $\nu\to 0$. This type of numerical setup still poses a significant challenge for the specialists in computational methods.

The principal problem consists in the fact that the analytical continuation of solution of the viscous BE, as a function of $x$, often has an infinite countable number of singular points which are simple poles with residues equal to $-2\nu$ that move in the complex plane while the time passes \cite{Four, Senouf1, Senouf2, BE_CP}. As the value of $\nu$ decreases, the singular points condense and their trajectories approach the segment of the real axis where the solution is sought. This leads to growth of the gradients of solution and of its derivatives. The example of such situation can be observed when studying the following exact solution of BE, see \cite{Dogan}:
\begin{equation}
\label{exSol}
u(t,x) = \frac{x/t}{1+\sqrt{t/t_0}\exp[x^2/(4\nu t)]},~~t\geq t_0>0, ~t_0=\exp\left(\dfrac{1}{8\nu}\right),
\end{equation}  
which has a countable number of simple poles:
\begin{equation}
\label{sing_exSol}
z_k^{\pm}=\pm2\sqrt{\nu t\big(\ln\sqrt{t_0/t}+i\pi (2k+1)\big)},~~k\in\mathbb{Z}.
\end{equation} 
It is easy to see that at $t=t_0$ the distance from the singular point $z_0^{\pm}$ to the real axis is equal to $\sqrt{2\pi\nu t_0}$. When $\nu$ is small, it leads to formation of a steep front of the solution. Similar effects appear for many of the known exact solutions of the BE, see the solutions from \cite{ExSols}. As time passes the singular points can either move away from the real axis (as in the above example) providing the reduction of the gradient, or approach it leading to an increase of the gradient.   

A systematic study of evolution of singular points in solutions of the BE started in \cite{SulemSulemFrish}, where the asymptotic analysis of the Fourier spectra of numerical solutions obtained at different time moments was performed. For the initial data $u(0,x)=-\sin(x)$ and $\nu=0.05$ it was shown that the solution develops a singularity of the type "pole" which, as $t$ increases, enters a small vicinity of the point $x=0$. It leads to formation of a steep gradient. However, in the case $\nu = 0$, as $t \to 1$, the singular points condense into a cubic root singularity and tend to $x=0$. This leads to the destruction of the classical solution, which corresponds to the breaking of the wave front (see for example Fig.~3 in \cite{LipKars}). To study the dynamics of singularities, besides analysing the Fourier spectra of the solution, one can apply Padé approximations. In this case, the zeros of their denominators can be used to approximate the locations of the singularities. Further developments of the mentioned methods and results can be found in \cite{FrishBec, Calfisch, Weidemann1, Weidemann2}.
An analytical approach for studying the dynamics of singularities of solutions to the BE in the case $\nu=0$ was developed in \cite{Zhuravleva}.

One should take into account that passing from the viscous to the inviscid case simply by setting $\nu\to 0$ raises a number of problems, whose treatment can be especially non-trivial when the spectral method with Galerkin truncation procedure is used (see \cite{Nazarenko1, Nazarenko2} and the references therein). In this context, it is very important to distinguish truncation-error effects from the real behaviour of the fluid as encoded in the equation. This means that we need to resolve very high-frequency oscillations accurately and with guaranteed precision.

From the mentioned facts it follows that the linear spaces of polynomials, in which solutions to the BE are usually sought, are not appropriate to construct a fast and accurate method. Indeed, following the direct and converse theorems of the approximation theory (see Bernstein's theorems, \cite{Bern} and chapter~8 in \cite{Tref}), one can conclude that for analytic solutions having singular points in a small vicinity of the real-valued interval of the problem, the convergence of polynomial approximations is formally exponential but so slow that one would need a mesh of many thousands of nodes to achieve adequate accuracy. For this reason, new methods of approximation and the algorithms for working with them are required. Since the analytic continuation of solution of the BE has simple poles, the natural choice is the application of rational approximations (RA). In this article we employ such approximations in combination with the collocation method. 

Note, however, that such combination can lead to a number of problems. The emergence of a pole of RA lying between the collocation nodes leads to large deviations of the solution and its approximation. Moreover, the algorithms based on the RA require an additional linearization with respect to the coefficients of RA. We should notice that the stable algorithms for solving differential equations using RA can be constructed by applying the least squares method \cite{Treff1}, which in combination with collocation method allows one to obtain a highly accurate solution with a small number of coefficients in the approximant (small degree of freedom)  \cite{ShapVP}. In this paper we employ another approach that enables one to cope with the mentioned difficulties. Based on it, we propose a new algorithm and test it by solving the BE for small $\nu$.   

Our algorithm is based on the application of rational barycentric interpolations (RBIs) \cite{Schneider, Baltensperger} and of the collocation method. The efficiency of this combination is ensured by the fact that for a given set of interpolation nodes the RBI is a linear combination of the values of the solution in these nodes and that the poles of RBI lay outside the interval of interpolation \cite{BerMit, FloaterHor}. It is worth noting that the RBIs have been already used for solving of a number of problems with singularities including the problems for the BE and its generalizations \cite{TeeTref, Idim, Bryndin, Sem22}. In these studies, for adapting the interpolations to singularities, special mappings were used that enabled one to move the singular point far away from the segment of the problem in the complex plane. For constructing such mappings, it is necessary to localize the singular point, which can be done by applying the Chebyshev--Pade approximations. However, the accuracy of such localization and its applicability when a singular point approaches extremely close to the domain of the solution were not studied yet.

An essential novelty of our work consists in the following three aspects.  First is the study of approximating properties of RBIs in cases where a singular point lies at a distance of the order $10^{-8}$--$10^{-4}$ from the domain of the problem. Such a study required the creation of a new numerical-analytical approach for computing the exact solutions of BE (see section~\ref{S2}) and the construction of a new mapping for adapting RBI to the singularity.

Second is a new algorithm based on the collocation and relaxation methods, which uses special modifications of polynomial and rational interpolations with respect to time and space, correspondingly. Hereinafter, we refer to it as the Spectral Time and Space Method or STSM. The accuracy and time performance of this method in the case of small $\nu$ were compared with those of a method based on the high-order Runge--Kutta scheme (RK5(4) 7S), which was also combined with RBIs \cite{Bryndin}. We characterize the cases where the STSM is the most effective method and those in which RK5(4) 7S exhibits higher performance.

Third, the use of the designed algorithms, for the first time, made it possible to numerically verify the accuracy of computing the trajectories of singular points for extremely small values of $\nu\sim 10^{-4}$--$10^{-6}$. This accuracy turns out to be rather high. Thus, the proposed methods represent an efficient numerical tool for studying the evolution of singular points in non-linear problems of hydrodynamics.

We close the introduction section by noting that the reasonable choice of one of the two mentioned methods enables one to achieve a better balance between the accuracy of numerical simulations and the computational costs in problems involving singularities. By this criterion, we compared the results obtained using the developed methods with those reported by other authors for two well-known test problems for the BE and ensured the principal advantage of the proposed methods.

\section{Statement of the problems and computation of their exact solutions}
\label{S2}

Before passing to development of the computational algorithm, we formulate the classical problems for the BE; give the formulas of their exact solutions; discuss the problems that arise when trying to use these formulas to perform computations in floating-point arithmetic; and propose a way to overcome them.  

The classical formulations for BE \eqref{BE} include the Cauchy problem (CP) \eqref{Cauch} posed for $x \in \mathbb{R}$ and the initial boundary-value problem (IBVP) \eqref{Cauch}, \eqref{BVP} posed in a certain segment (hereinafter, we set $x\in[-1,1]$):
\begin{equation}
\label{Cauch}
u(0,x) = \varphi(x),
\end{equation}
\begin{equation}
\label{BVP}
u(t,-1) = 0 , \quad u(t,1) = 0.
\end{equation}
In what follows we shall assume that $\varphi(x)$ is at least continuous. In our numerical tests we consider the case of analytic initial data.

The solution to \eqref{BE} can be expressed using the Cole--Hopf transformation  \cite{Hopf, Cole}:
$$u(t,x) = -2 \nu \dfrac{\partial}{\partial x}\ln{F(t,x)} = -2 \nu \dfrac{F_x}{F},$$
which reduces the CP \eqref{BE}, \eqref{Cauch} to the Cauchy problem for the heat conduction equation \eqref{HE}, \eqref{CauchHE} and transforms the IBVP \eqref{BE}, \eqref{Cauch}, \eqref{BVP} into the problem \eqref{HE}--\eqref{BVPHE}:  
\begin{equation}
\label{HE}
 F_t - \nu F_{xx} = 0,\quad   t \in [T_0, T],
\end{equation} 
\begin{equation}
\label{CauchHE}
F(0,x) = F_0(x) = \exp \left[-\dfrac{1}{2\nu} \int_0^x \varphi(\xi) d\xi \right],
\end{equation} 
\begin{equation}
\label{BVPHE}
F_x(t,-1) = 0 , \quad F_x(t,1) = 0.
\end{equation}

The solution to the CP \eqref{HE}, \eqref{CauchHE} on the whole real axis can be written using the Poisson integral
\begin{equation}
\label{eq_tepl_answer}
	F(t,x) = \dfrac{1}{\sqrt{4\nu\pi t}} \int_{-\infty}^{\infty} \exp \left[ -\dfrac{(x-y)^2}{4\nu t} -\dfrac{1}{2\nu} \int_0^y \varphi(\xi) d\xi \right] dy.
\end{equation}

The solution to the IBVP \eqref{HE}--\eqref{BVPHE} can be obtained by the Fourier method:
\begin{equation}
\label{FourierSol}
F(t,x) = A_0+2\sum\limits_{k=1}^{\infty}A_k\cos(k\pi x)\exp(-k^2\pi^2\nu t),~~~A_k=\int\limits_0^1\cos(k\pi x)F_0(x)dx.
\end{equation}

Hence, the solutions to problems \eqref{BE}, \eqref{Cauch} and \eqref{BE}, \eqref{Cauch}, \eqref{BVP} can be expressed by the following formulas, respectively:  
\begin{equation}
\label{eq_ub_exact}
	u_{ex}(t,x) = \dfrac
{\int_{-\infty}^{\infty}\dfrac{x-y}{t} \exp \left[ -\dfrac{(x-y)^2}{4\nu t} -\dfrac{1}{2\nu} \int_0^y \varphi(\xi) d\xi  \right] dy}
{\int_{-\infty}^{\infty} \exp \left[ -\dfrac{(x-y)^2}{4\nu t} -\dfrac{1}{2\nu} \int_0^y \varphi(\xi) d\xi \right] dy},
\end{equation} 
\begin{equation}
\label{eq_ub_exact_Four}
	u_{ex}(t,x)  = \dfrac
{4\pi\nu\sum\limits_{k=1}^{\infty}kA_k\sin(k\pi x)\exp(-k^2\pi^2\nu t)}
{A_0+2\sum\limits_{k=1}^{\infty}A_k\cos(k\pi x)\exp(-k^2\pi^2\nu t)}.
\end{equation}

\begin{remark}
\label{rem1}
The CP \eqref{Cauch} and the IBVP \eqref{Cauch}, \eqref{BVP} have a completely different physical sense. Nevertheless, for testing numerical methods, it is useful to note that if the initial condition is such that the functions
\begin{equation}
\label{PhiFun}
\tilde{\Phi}_{\pm}(s)=\int\limits_{0}^{\pm 1-s}\varphi(\xi)d\xi~~\text{are even functions}
\end{equation}
of the variable $s$, then the solution \eqref{eq_ub_exact} satisfies the boundary conditions \eqref{BVP}. In such a case the solutions \eqref{eq_ub_exact} and \eqref{eq_ub_exact_Four} coincide on the segment $x\in[-1,1]$.
\end{remark}
\begin{remark}
\label{rem2}
It is worth noting that, to keep the remainders of the series in \eqref{eq_ub_exact_Four} small, one should take
$k>1/(\pi\sqrt{\nu t})$. Hence, to obtain a highly-accurate result for sufficiently small values of $t$ and $\nu$, one may need to compute a large number of terms. In this case, the use of the expression for the coefficients $A_k$ from \eqref{FourierSol} leads to a problem of computing the integral of the rapidly oscillating function. In our opinion, an efficient way for searching $A_k$ in this case is the application of the Fast Fourier Transform to the function $F_0(x)$.
\end{remark}

Notice that despite a considerable simplicity of formulas \eqref{eq_ub_exact}, \eqref{eq_ub_exact_Four}, their application in the case of small $\nu$ and of initial data $\varphi(x)$ having negative values in a vicinity of the point $x=0$ causes a problem. Indeed, assume that there exists such point $y\in I$ (here, $I$ is the domain of problem) that 
$$\Phi(y)=\int_0^y \varphi(\xi)d\xi<-\Phi_M,$$
where $\Phi_{M}=600\nu\ln 10$. As a result, the quantities of order higher than $10^{300}$ emerge in formulas \eqref{CauchHE}, \eqref{eq_tepl_answer}, \eqref{eq_ub_exact}. In floating-point arithmetic this causes an overflow.
The reason is that for $y=x$ the argument of exponential functions in \eqref{eq_tepl_answer}, \eqref{eq_ub_exact} exceeds the quantity $\Phi_M/2\nu$, hence, the values of exponents exceed $10^{300}$. The same is valid for the exponent standing in \eqref{CauchHE}. This problem arises when, for example, $\varphi(x)=-\sin \pi x$ and $\nu\leq0.4\times 10^{-3}$. An initial data of this form was, in particular, considered in \cite{Basdevan,OVasilyev,SulemSulemFrish}. In \cite{Basdevan,OVasilyev} for searching the exact solutions the formulas \eqref{eq_ub_exact}, \eqref{eq_ub_exact_Four} were used. However, the comments related to the mentioned problem are absent.

The problem can also arise for positive initial data: if $\exists x\in I: \forall y\in \mathbb{R}$ $(x-y)^2/(2t)+\Phi(y)>\Phi_M$, then in floating-point arithmetic the numerical values of the numerator and denominator of \eqref{eq_ub_exact} at such a point $x$ turn out to be zero. 

Note also that the variations of the values of function by quantities of the order $10^{15}$ and higher, when $x$ varies over $I$, are undesirable in double-precision format. The reason is that applying such operations as differentiation, integration, and many others to such functions requires computing linear combinations of their values at different nodes of $I$. Since the mantissa of a real number in double precision is limited to 16 significant digits, operations of this kind can lead to catastrophic growth of round-off error. For achieving high accuracy, it is important that the values of functions standing in the numerators and denominators of \eqref{eq_ub_exact}, \eqref{eq_ub_exact_Four} remain bounded in the order of magnitude as $x$ varies over $I$.

Taking into account the above arguments, in order we can use \eqref{eq_ub_exact}, let us impose the constraints on the initial data:
\begin{equation}
\label{bound}
\forall t>0, y\in I~~\Phi(y)>-\Phi_m \text{~~and~~} \forall t>0, x\in I~~\exists y\in \mathbb{R}: (x-y)^2/(2t)+\Phi(y)<\Phi_m,
\end{equation}
where $\Phi_m=15\nu \ln 10$. While writing the last inequality we assume that the initial data can be continued from the domain $I$ into $\mathbb{R}$. Let us introduce the function
$$
B(t,x)=\min\limits_{y\in\mathbb{R}}\bigg(\frac{(x-y)^2}{2t}+\Phi(y)\bigg).
$$
It is worth noting that for $t>0$, the function $\frac{(x-y)^2}{2t}+\Phi(y)$ should tend to infinity sufficiently fast when $|x|\to\infty$. This constitutes the conditions of existence of the solution to BE. Together with the requirement that $\Phi(y)$ is continuous (which follows directly from the continuity of $\varphi(x)$), this implies that the minimum appearing in the formula for $B(t,x)$ exists. Therefore, when $y$ varies from $-\infty$ to $\infty$ the function $\frac{(x-y)^2}{2t}+\Phi(y)$ decreases from infinity, can oscillate about some level but attaining a global minimum, and then goes back to infinity.

Then, due to the fact that the term $(x-y)^2/2t$ is non-negative and is zero only for $x=y$, the first condition in \eqref{bound} can be rewritten as $\forall t>0, x\in I~~ B(t,x)>-\Phi_m$. The second condition can be reformulated as $\forall t>0, x\in I~~B(t,x)<\Phi_m$, and instead of \eqref{bound} we can write
\begin{equation}
\label{bound_2}
\forall t>0, x\in I~~|B(t,x)|<\Phi_m.
\end{equation}

If the values of $\nu$ and $\varphi(x)$ are such that one of the conditions \eqref{bound} fails then we propose to multiply \eqref{eq_ub_exact} by the function $\exp[\psi(t,x)/(2\nu)]$, where $\psi(t,x)$ should be chosen so that for the function
$\Psi(t,x)=B(t,x)-\psi(t,x)$
the following condition holds: 
\begin{equation}
\label{Fbound}
\forall t>0, x\in I~~|\Psi(t,x)|<\Phi_m.
\end{equation}

For the most of initial data, it is impossible to express the function $B(t,x)$ analytically, therefore determining $\psi(t,x)$ can be based on symbolic computations which are free of the restrictions inherent in floating-point arithmetic. Let us describe the way of computing $\psi(t,x)$ on the example of problem with the initial data $\varphi(x)=-\sin(\pi x)$, $x\in[-1,1]$. For simplicity we fix $\nu=10^{-4}$, $t=1$. In this case, the function $\psi(x):=\psi(1,x)$ can be expressed as
\begin{equation}
\label{psi1}
\psi(x)=-(2\nu M\ln 10)\overline{\psi}_P(x),
%0.375128(x-\sign(x))^2-0.63581.
\end{equation}
where the quantity $M$ is described below,
$$\overline{\psi}_P(x)=-0.59(x-\sign(x))^2+1.$$
Let us explain, how this expression was obtained. Denote the common logarithm of the denominator of \eqref{eq_ub_exact} by $D_p(x)$; $\overline D_p(x) = D_p(x)/M$, where $M = \|D_p(x)\| = D_p(1)$, $\|\cdot\|$ is the supremum norm of the continuous function in the segment $[-1,1]$. 

\begin{figure}[h!]
\centering
\includegraphics[width=1.0\linewidth]{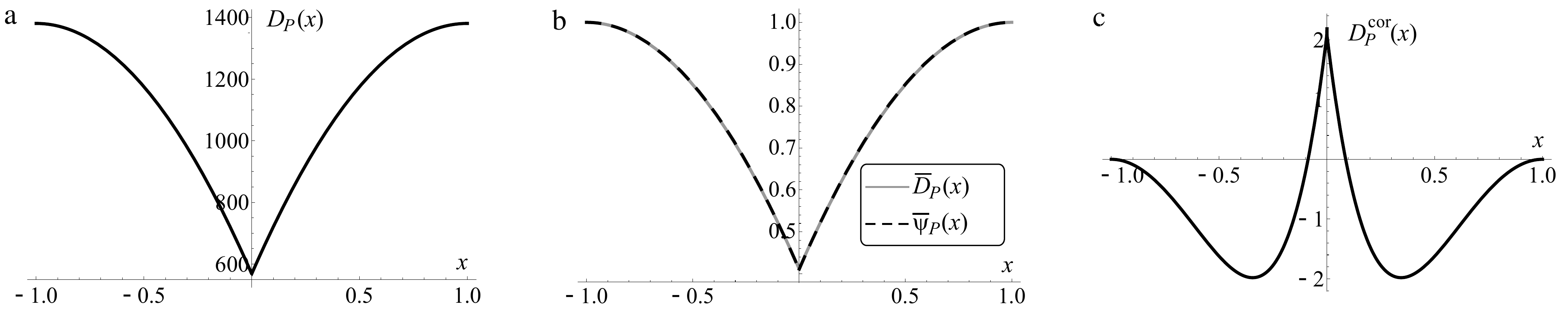}
\caption{Denominator of the solution \eqref{eq_ub_exact} for $\varphi(x)=-\sin \pi x$ and $\nu=10^{-4}$ and the correction of its values: plots of the function $D_P(x)$ (a); the function $\overline{D}_P(x)$ and its approximation $\overline{\psi}_P(x)$ (b); and the function $D_P^{\mathrm{cor}}(x)$ obtained using the correcting function (c).}
\label{fig1_0}
\end{figure}

Denote by
$$D_p^{\mathrm{cor}}(x) = D_p(x)+\psi(x)/(2\nu\ln 10)$$
the common logarithm of the product of denominator of \eqref{eq_ub_exact} and the correcting function $\exp(\psi(x)/(2\nu))$ 
The plot of the function $D_p(x)$ is given in Fig.~\ref{fig1_0},~a. It can be seen that the values of the denominator of \eqref{eq_ub_exact} exceed the quantity $10^{300}$ by many orders of magnitude and the variation of these values as $x$ ranges over $I$ exceeds 700 orders. Such situation is completely unsuitable for performing the computations in floating point systems. 

The approximation of the normalized function $\overline D_p(x)$ by the polynomial $\overline \psi_P(x)$ (see Fig.~\ref{fig1_0},~b) enables one to compute $\psi(x)$.
Indeed, let the relative error of the mentioned approximation be 
$$\varepsilon_P=\|\overline D_p(x) - \overline{\psi}_P(x)\|/\|\overline D_p(x)\|\ll 1.$$ 
Hence, if similarly to \eqref{psi1}, one takes $\psi(x)=-(2\nu M\ln 10)\overline{\psi}_P(x)$, then
$$
\begin{aligned}
&\|D_p^{\mathrm{cor}}(x)\|=\|D_P(x)-M\overline{\psi}_P(x)\|=M\|\overline D_p(x) - \overline{\psi}_P(x)\|\leq \\
&\leq \varepsilon_P M\max\limits_{x\in[-1,1]}|\overline D_p(x)|=\varepsilon_P\max\limits_{x\in[-1,1]}|M\overline D_p(x)|=\varepsilon_P\|D_p(x)\|\ll \|D_p(x)\|.
\end{aligned}
$$
The plot of $D_p^{\mathrm{cor}}(x)$ is given in Fig.~\ref{fig1_0},~c.

\begin{figure}[h!]
\centering
\includegraphics[width=0.8\linewidth]{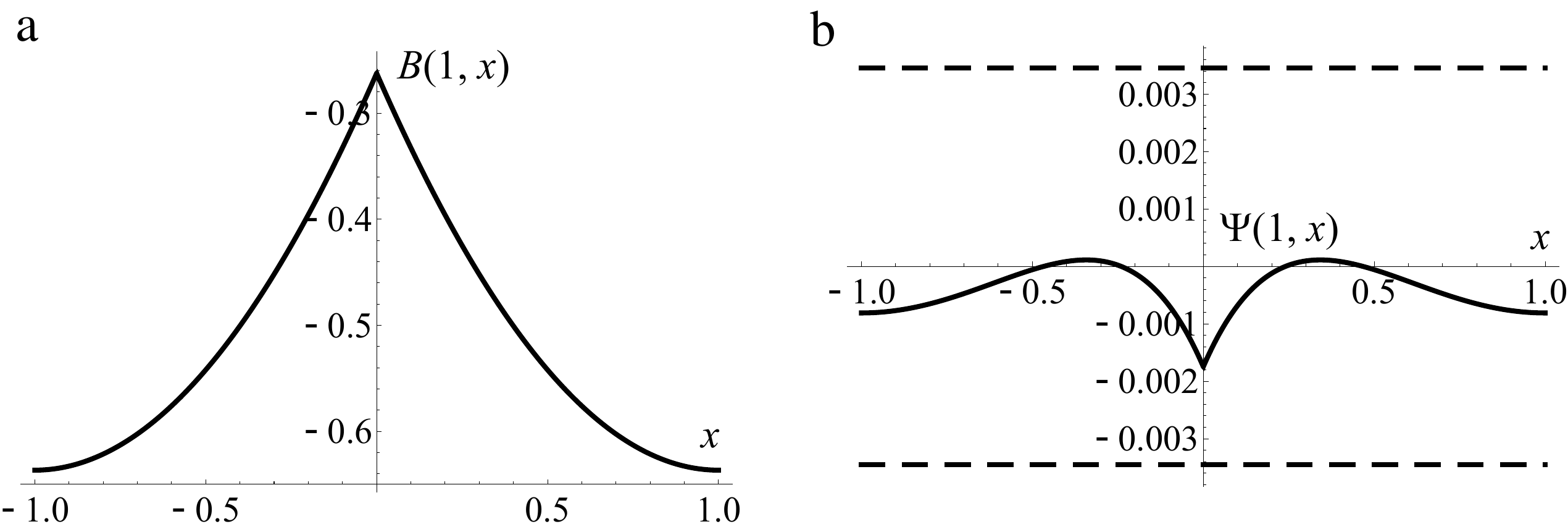}
\caption{Correction of the arguments of exponential functions in \eqref{eq_ub_exact} for $\varphi(x)=-\sin \pi x$ and $\nu=10^{-4}$: plots of the function $B(1,x)$ (a) and the function $\Psi(1, x)$ -- the restrictions $\pm\Phi_m$ of inequality \eqref{Fbound} are shown by dashed lines (b).}
\label{fig1_1}
\end{figure}

Note that the use of the correcting factor $\exp(\psi(x)/(2\nu))$ ensures that \eqref{Fbound} holds. Figure \ref{fig1_1} contains the plots of functions $B(1,x)$ and $\Psi(1,x)$. One can see that restrictions \eqref{bound_2} fail for the function $B(1,x)$, whereas for the function $\Psi(1,x)$ all the restrictions \eqref{Fbound} are satisfied (in this example $\Phi_m\approx 0.0035$). Moreover, the correcting factor balances the values of the numerator and denominator of \eqref{eq_ub_exact} in such a way that the range of their orders of magnitude is much less than 15 as $x$ varies over $I$: here the range is approximately $10^{-2}$ -- $10^2$, see Fig.~\ref{fig1_0},~c, i.e. it is equal to 4 orders of magnitude.

In the general case, where $\varphi(x)=-\sin \pi x$, $t>0$ and  $\nu$ is arbitrary small, we shall use the following series for approximating the function $\overline D_p(x)$:  
    $$
\overline{\psi}_P(t,x)=\sum\limits_{k=0}^{\overline{K}}A_{\psi,k}T_k(2|x|-1),
    $$
    where $T_k$ is the Chebyshev polynomial of the first kind of the power $k$. By writing the approximation in this form, we have taken into account that the function appearing in the denominator of \eqref{eq_ub_exact} is even (it is based on the fact that $\Phi(y)$ is even). Further, we set again $\psi(x)=-(2\nu M\ln 10)\overline{\psi}_P(x)$. Numerical experiments show that truncating the series at such $k=\overline{K}$ that $|A_{\psi,k}|<10\nu$ suffices to satisfy \eqref{Fbound}. To compute the coefficient $A_{\psi,k}$ we developed the following Wolfram Mathematica code: 

\begin{footnotesize}
\begin{lstlisting}
nu := 10^(-5);
t := 0.5;
Nn0 := 3;
coefs = ConstantArray[1, Nn0-2];
Nn = Nn0;
DP[x_]:=Quiet[NIntegrate[Exp[-(x-y)^2/(4*nu*t)-(Cos[Pi*y]-1)/(2*nu*Pi)], 
{y, -Infinity, Infinity}, Method -> "ExtrapolatingOscillatory", 
WorkingPrecision -> 35, PrecisionGoal -> 15, MaxRecursion -> 20]];
M := Log10[DP[1]];
DDP[x_] := Log10[DP[x]]/M;
While[Abs[Last[coefs]] > nu*10,
   Nn = Nn+2;
   NodesF[j_] := (Cos[(2*j-1)*Pi/(2*Nn)]+1)/2;
   Nodes := Array[NodesF, Nn, {1, Nn}];
   BaseF[j_] := ChebyshevT[j, 2*x-1];
   Base := Array[BaseF, Nn, {0, Nn-1}];
   DDPL := Transpose[{Nodes, DDP[Nodes]}];
   FFcor[x_] = Fit[DDPL, Base, x];
   coeF[j_] := Coefficient[FFcor[x], ChebyshevT[j, 2*x-1], 1];
   coefs = Array[coeF, Nn-2, {2, Nn-1}];]
ress[x_]:=FFcor[x]-Sum[coefs[[j-1]]*ChebyshevT[j,2*x-1], {j,2,Nn-1}];
coefs = Join[{Coefficient[ress[x], ChebyshevT[1, -1+2*x], 0], 
    Coefficient[ress[x], ChebyshevT[1, -1+2*x], 1]}, coefs];
Apsi = -2*nu*M*Log[10]*coefs;
Print["Number of coefficients of expantion over ", T_k, " is overline{K}+1=", Nn]
Print["The coefficients are \n", Apsi]

(*You can save the coefficients to file*)
(*Export["c:/coefs.txt", Apsi, "TSV"];*)

(*You can plot the error of approximation of \overline{D}_P by \overline{\psi}(t,x).
This can take a time*)
(*Plot[DDP[x]-FFcor[x], {x,0,1}, PlotStyle->{Thickness[0.007`],Black},
LabelStyle -> {FontFamily -> "Times New Roman", 20, GrayLevel[0]}, 
ImageSize -> Large, PlotLegends -> Placed[LineLegend[{"Error of
approximation of DP by psi(t,x)"}, LegendFunction -> Frame, 
LegendMarkerSize -> {{55, 10}}], {0.5, 0.2}]]*)

(*You can plot the value of \Psi(t,x). This can take a time*)
(*FFcor2[x_] := Sum[Apsi[[j]]*ChebyshevT[j-1, 2*Abs[x]-1], {j, 1, Nn}];
Psi[x_] :=  MaxValue[-(x-y)^2/(4*nu*t)-(Cos[Pi*y]-1)/(2*nu*Pi)+FFcor2[x]/(2*nu), y];
Plot[Log10[Abs[Psi[x]]],{x,-1,1},PlotStyle->{Thickness[0.007`],Black},
LabelStyle -> {FontFamily -> "Times New Roman", 28, GrayLevel[0]},
ImageSize -> Large, PlotPoints -> 200]*)
\end{lstlisting}
\end{footnotesize}

\indent In this work we develop a numerical method for solving the IBVP for the BE:
\begin{equation}
\label{eq_ub} 
\begin{cases}
	u_t + u u_x = \nu u_{xx}, \quad x \in [-1,1] ,
	\\
	u(0,x) = \varphi(x),
	\\
	u(t,-1) = 0 , \quad u(t,1) = 0,
\end{cases}
\end{equation}
with the following variants of setting the initial data
\begin{eqnarray}
\label{eq_ub1}
\varphi(x) &=& -\sin\pi x,\quad    t \in [0,T] ~~\text{ (here $T_0=0$)};\\
\label{eq_ub2}
\varphi(x) &=& \frac{x+1}{1+ \exp[(x+1)^2/(4\nu) -  1/(16\nu) ]}, \quad   t \in [1,T]~~\text{ (here $T_0=1$)}.
\end{eqnarray} 
Notice that the initial data \eqref{eq_ub2} corresponds to the exact solution \eqref{exSol} with the shift of the variable $x \in [0,2] \Rightarrow x \in [-1,1]$.

\begin{figure}[h!]
\centering
\includegraphics[width=0.80\linewidth]{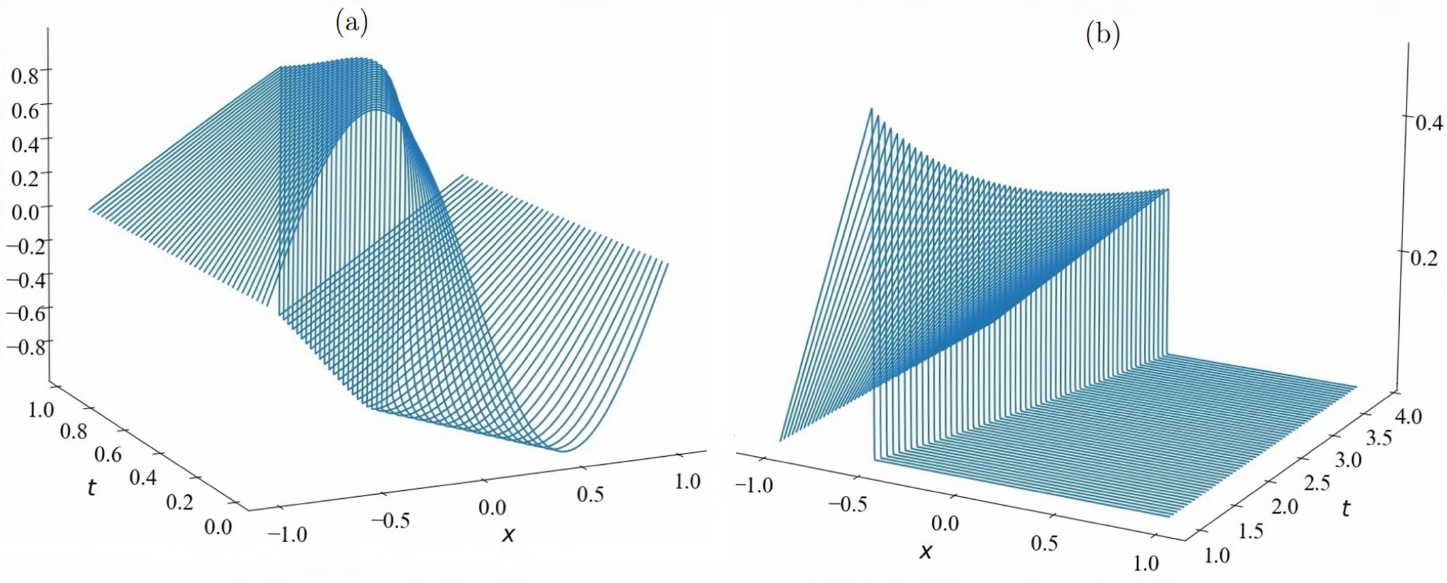}
\caption{Plots showing the evolution of numerical solutions to \eqref{eq_ub} as $\nu=10^{-4}$ for the initial data \eqref{eq_ub1} (a) and for the initial data \eqref{eq_ub2} (b).}
\label{fig1}
\end{figure}

\indent In Fig.~\ref{fig1} we put the plots showing the evolution of numerical solutions to the problem \eqref{eq_ub} with the initial data specified above, which were obtained using the algorithm described below. Note that, due to remark~\ref{rem1} the exact solution of IBVP \eqref{eq_ub}, \eqref{eq_ub1} can be expressed by formula \eqref{eq_ub_exact}.

\section{Algorithm based on barycentric interpolation}
\label{S3}

\indent In this section we describe the interpolation formulas for approximation of the solution to the problem \eqref{eq_ub} and give the expressions for elements of the matrices approximating the differentiation operators with respect to the variables $x$ and $t$. To find the solution, we develop a numerical algorithm based on the relaxation technique and on the application of the collocation method. The idea of the latter consists in substituting the approximation into the equation of the problem and in writing the obtained equation at the collocation points.  It is worth noting that the approximation with respect to the time variable takes the initial condition into account directly.

To approximate the solution with respect to the variable $x \in [-1,1]$ for a fixed $t$, we can start with the rational barycentric interpolation (RBI)
\begin{equation}
\label{bar0}
	u(t,x)\approx R_N(x) =  \sum\limits_{j=1}^N \dfrac{\omega_j u(t, g(x_j^{\mathrm{ch}}))}{x-g(x_j^{\mathrm{ch}})} \bigg/ \sum\limits_{j=1}^N\dfrac{\omega_j}{x-g(x_j^{\mathrm{ch}})},
\end{equation}
where $x_j^{\mathrm{ch}}=\cos\dfrac{(2j-1)\pi}{2N}$ and $\omega_j = (-1)^{j-1}\sin\dfrac{(2j-1)\pi}{2N}$ are the nodes and the weights of the interpolation, respectively, $j = \overline{1,N}$, $g(x)$ is the conformal mapping that maps the segment $[-1,1]$ onto itself and allows one to adapt the nodes of interpolation to the location of the singular point, $g(\pm 1)=\pm 1$. In what follows, we set $x_j=g(x_j^{\mathrm{ch}})$. Note that the approximation \eqref{bar0} was proposed in \cite{Schneider}. The error estimates for it were obtained in \cite{Baltensperger} and in later papers of professor J.-P. Berrut and his co-authors. The efficiency of application of \eqref{bar0} for solving BE was demonstrated in \cite{Idim, TeeTref}.

\indent We shall seek the solution of the problem \eqref{eq_ub} in the computational cell $C_k=\{(t,x): t \in [\tau_k, \tau_{k+1}], x \in[-1,1] \}$, where the segment $ [\tau_k, \tau_{k+1}]$ will be described later, in the form of tensor product of the following interpolations: 

\begin{itemize} 
\item Modified interpolation \eqref{bar0} for a fixed $t \in [\tau_k, \tau_{k+1}]$:
\begin{equation}
\label{eq_baricentr}
	\widehat{R}_N(x) =  \sum\limits_{j=1}^N  \dfrac{1-x^2}{1-x_j^2}  ~\dfrac{\omega_j u(t, x_j)}{x-x_j} \bigg/ \sum\limits_{j=1}^N\dfrac{\omega_j}{x-x_j},~~x_j=g(x_j^{\mathrm{ch}}).
\end{equation}
\item Modified interpolation polynomial with the Chebyshev nodes, $\theta_i=\cos\dfrac{(2i-1)\pi}{2M}$, and the weights $\xi_i = (-1)^{i-1}\sin\dfrac{(2i-1)\pi}{2M}$, $i = \overline{1,M}$, for $t \in [\tau_k, \tau_{k+1}]$ with a fixed $ x \in[-1,1]$:
\begin{equation}
\label{eq_intetpolate}
	\widehat{P}_M(t) = \dfrac{1}{M}\sum\limits_{i=1}^M \dfrac{1+\Lambda(t)}{1+\theta_i} \dfrac{\xi_i T_M(\Lambda(t))}{\Lambda(t) - \theta_i}u(t_i, x) + u(\tau_k,x), 
\end{equation}
   where $\Lambda(t)$ is a linear mapping of the segment $[\tau_k, \tau_{k+1}]$ onto the segment $[-1,1]$, $t_i =\Lambda^{-1}(\theta_i)$, and $T_M$ is a Chebyshev polynomial of degree $M$. 
\end{itemize}

In the case of \eqref{eq_baricentr} the modification of the interpolation consists in the presence of the multiplier $\dfrac{1-x^2}{1-x_j^2}$, $x_j=g(x_j^{\mathrm{ch}})$, which allows one to take into account zero boundary conditions, whereas in the case of \eqref{eq_intetpolate}, it consists in the presence of the multiplier $\dfrac{1+\Lambda(t)}{1+\theta_i}$ and the term $u(\tau_k,x)$, which ensures the realization of the initial condition. It should be noted that there are other techniques that allow one to account for the boundary conditions and initial data when spectral methods are used. For a discussion of some of the ideas, we refer the reader to \cite{Allal}.

\begin{remark}
\label{rem3}
Introducing the mentioned multipliers plays a fundamental role in developing the algorithm based on the collocation method. Indeed, direct application of \eqref{bar0} leads to a problem of linear algebra, where the matrices approximating the derivatives (see \eqref{TeeTreff1}, matrices $D_{(1)}$ and $D_{(2)}$) have only zero eigenvalues and the corresponding Jordan blocks of large sizes. This is a direct consequence of the fact that the differentiation operators of the first and second orders in the space of polynomials have only constant and linear eigenfunctions, which correspond to zero eigenvalues. Therefore, direct usage of the matrices $D_{(1)}$ and $D_{(2)}$ is completely unsuitable for organizing computations.
\end{remark}

\begin{remark}
If $g(x)$ is the identity mapping, then the approximation \eqref{eq_baricentr} is a polynomial; otherwise \eqref{eq_baricentr} is a rational function (RF). In what follows, we describe an algorithm for adapting the spatial grid to a singular point. This algorithm employs a special form of the mapping $g(x)$. 
\end{remark}

Due to the use of \eqref{eq_intetpolate} in each cell $C_k$, the time grid is defined by the Chebyshev nodes $\theta_i$, $i = \overline{1,M}$, which are clustered near the boundaries, but none of them fall on the boundary exactly. The proposed numerical algorithm uses the collocation method and, hence, it allows one to find approximate values of the solution at these nodes. However, to set the initial data for computing the solution in the cell $C_{k+1}$, we need to calculate the values of the solution at the right boundary of the considered time segment, i.e. at $t=\tau_{k+1}$. To this end, we use the interpolation polynomial with Chebyshev nodes without modifications:
\begin{equation}
\label{eq_int}
	u(\tau_{k+1},x) \approx P_M(\tau_{k+1}) = \dfrac{1}{M}\sum\limits_{i=1}^M  \dfrac{\xi_i}{1 - \theta_i}u(t_i,x), 
\end{equation} 
where $x\in[-1,1]$ is fixed.

For getting the matrices approximating the differentiation operators with respect to the variables $t$ and $x$, we need to find the first and the second derivatives of the interpolation formulas, and then take the limits $t \rightarrow t_i $ and $x \rightarrow x_j$. The elements of the matrices $D_{(1)}$ and $D_{(2)}$, which approximate the first and second order derivatives, can be derived using \eqref{bar0} (see \cite{TeeTref} for instance):
\begin{equation}
\label{TeeTreff1}
	D_{(1) ij} =
\begin{cases}
	\dfrac{\omega_j}{\omega_i(x_i-x_j)}, &\text{at } i \neq j;
	\\
    - \chi_i, &\text{at } i = j;
\end{cases}~~~
	D_{(2) ij} =
\begin{cases}
	\gamma_{ij}, &\text{at } i \neq j;
	\\
    - \sum\limits_{k=1,k\neq i}^N \gamma_{ik}, &\text{at } i = j;
\end{cases}
\end{equation}
where $ \chi_i = \sum\limits_{k=1,k\neq i}^N \dfrac{\omega_k}{\omega_i(x_i-x_k)}$, $\gamma_{ik} = -\dfrac{2\omega_k}{\omega_i(x_i-x_k)} \left(  \chi_i  +\dfrac{1}{x_i-x_k}\right)$.

Now, let us write the matrix approximations of differentiation operators, which we obtain from the modified interpolation formulas:

\begin{itemize} 
\item $B_{(1)}=(B_{(1) ij})$ is an $N\times N$ matrix approximating the first-order derivative with respect to $x$, obtained from \eqref{eq_baricentr}: 
\begin{equation}
\label{eq_diffX1}
	B_{(1) ij} =
\begin{cases}
	\dfrac{1-x_i^2}{1-x_j^2}\dfrac{\omega_j}{\omega_i(x_i-x_j)}, &\text{if}~~  i \neq j;
	\\
	-\dfrac{2x_i}{1-x_i^2} - \chi_i, &\text{if}~~  i = j.
\end{cases}
\end{equation}

\item $B_{(2)}=(B_{(2) ij})$ is an $N\times N$ matrix approximating the second-order derivative with respect to $x$, obtained from \eqref{eq_baricentr}:
\begin{equation}
\label{eq_diffX2}
	B_{(2) ij} =
\begin{cases}
	-\dfrac{4x_i}{1-x_j^2}\dfrac{\omega_j}{\omega_i(x_i-x_j)}+\dfrac{1-x_i^2}{1-x_j^2}\gamma_{ij},  &\text{if}~~ i \neq j;
	\\
	-\dfrac{2}{1-x_i^2} + \dfrac{4x_i}{1-x_i^2}\chi_i - \sum\limits_{k=1,k\neq i}^N \gamma_{ik}, &\text{if}~~  i = j.
\end{cases}
\end{equation}

\item $A= (A_{ij})$ is an $M\times M$ matrix approximating the first-order derivative with respect to $t$, obtained from \eqref{eq_intetpolate}:
\begin{equation}
\label{eq_diffT}
A_{ij} =
\begin{cases}
(-1)^{i+j}\dfrac{2}{\tau_{k+1} - \tau_k}\dfrac{(1+\theta_i)\sqrt{1-\theta_j^2}}{(1+\theta_j)(\theta_i-\theta_j)\sqrt{1-\theta_i^2}}, &\text{if}~~  i\neq j;
\\
\dfrac{2}{\tau_{k+1} - \tau_k} \left(\dfrac{1}{1+\theta_i}+\dfrac{\theta_i}{2(1-\theta_i^2)} \right), &\text{if}~~  i=j.
\end{cases}
\end{equation} 
\end{itemize} 

Before constructing the described approximations, we partitioned the domain $\hat \Omega = [T_0, T] \times [-1,1] $ of problem \eqref{eq_ub} into cells as follows: we introduce a uniform time grid with the steps $d_k$ and the nodes $\tau_k$ such that $\tau_0 = T_0$, $\tau_{k+1} = \tau_k + d_k$, $k = \overline{0,K}$ and $\tau_K = T$.

We find the solution in each cell $C_k =  [\tau_k, \tau_{k+1}] \times [-1,1]$ using the relaxation method and combining the collocation method with the tensor-product approximation of the solution constructed from \eqref{eq_baricentr} and \eqref{eq_intetpolate}. To this end, we introduce the grid $(t_i,x_j)$ in $C_k$, where $t_i =\Lambda^{-1}(\theta_i)$, $\theta_i=\cos\dfrac{(2i-1)\pi}{2M}$, $i=\overline{1,M}$. In the cell $C_0$ we set $x_j=x_j^{\mathrm{ch}}=\cos\dfrac{(2j-1)\pi}{2N}$, $j=\overline{1,N}$, while in other cells we set $x_j=g(x_j^{\mathrm{ch}})$, where $g(x)$ serves to adapt the grid to a singularity of the solution. 

To implement the relaxation method for computing the solution in the cell $C_k$, we introduce a fictive variable $s$ and set $u=u_0(t,x,s)+u_{C_k}(x)$, where $u_0$ satisfies a homogeneous initial condition at the point $t=\tau_k$ and $u_{C_k}(x)$ satisfies a non-homogeneous initial condition in the cell $C_k$. The values of $u_{C_k}(x)$ can be computed using the values of the solution obtained in the previous cell in accordance with \eqref{eq_int}. The relaxation method consists in finding the values of $u_0(t,x,s)$ as the limit of the solution to the problem:
\begin{equation}
	\dfrac{\partial u_0}{\partial s} =  u_{0t} + (u_0+u_{C_k}) (u_0+u_{C_k})_x - \nu (u_0+u_{C_k})_{xx} \text{ as } s\rightarrow\infty.
\end{equation} 
In what follows, we set $u_0(t,x,0)=0$ and introduce a uniform grid over $s$ with the step $\Delta$ and the nodes $s^n=n\Delta$, $n=0,1,2,...$. Next, we define $u_0^n=u_0(t,x,s^n)$, approximate $\dfrac{\partial u_0}{\partial s}\approx \dfrac{u_0^{n+1}-u_0^{n}}{\Delta}$, and derive the following iterative scheme of the relaxation method for searching the solution in cell $C_k$:
\begin{equation}
\label{eq_MetUst}
	\dfrac{u_0^{n+1}-u_0^{n}}{\Delta} = u_{0 t}^{n+1} + (u_0^n+u_{C_k}) (u_0^n+u_{C_k})_x - \nu (u_0^{n+1}+u_{C_k})_{xx}.
\end{equation}
We are now in position to apply the collocation method. Let 
\begin{itemize}

\item[$\circ$] $U_0^{n}$ be an $M \times N$ matrix containing the approximate values of the solution at the collocation nodes at the $n$th iteration of the relaxation method, with $M$ and $N$ being the numbers of nodes in time and space, respectively.

\item[$\circ$] $ U_{C_k}$ be an $M \times N$ matrix whose rows are identical and contain the values of the solution at the nodes $x_j$, $j=\overline{1,N}$, at $t=\tau_k$, obtained using the values of the solution from the previous cell.
\end{itemize}
We have the following approximate equalities:
\begin{itemize}
\item[$\circ$]  $(U_0^n+U_{C_k})_x \approx (U_0^{n} + U_{C_k})B_{(1)}^{T}$, where the elements of $B_{(1)}$ are defined in \eqref{eq_diffX1}.

\item[$\circ$] $(U_0^{n+1}+U_{C_k})_{xx} \approx (U_0^{n+1} + U_{C_k}) B_{(2)}^{T}$, where the elements of $B_{(2)}$ are defined in \eqref{eq_diffX2}.

\item[$\circ$] $ U_{0 t}^{n+1} \approx  AU_0^{n+1}$, where the elements of $A$ are defined in \eqref{eq_diffT}.
\end{itemize}
The subscripts ``$x$'', ``$xx$'' and ``$t$'' indicate that the matrices contain the values of the corresponding derivatives at the collocation nodes.

Writing \eqref{eq_MetUst} at the collocation nodes, we get the matrix equation:
\begin{equation}
\label{eq_Iter}
    U_0^{n+1}-\Delta  A U_0^{n+1}+  \Delta \nu U_0^{n+1} B_{(2)}^T = 
 U_0^{n}+\Delta \left[ (U_0^{n} + U_{C_k})\odot \big((U_0^{n} + U_{C_k})B_{(1)}^T\big)   -   \nu U_{C_k}B_{(2)}^T\right],
\end{equation}  
where the symbol "$\odot$" means element-wise multiplication of matrices.

Next, we use the spectral decomposition of the matrices $A$ and $B^T_{(2)}$:
\begin{equation}
\label{eq_spectral}
A = R_{A}D_{A}R_{A}^{-1},~~ B^T_{(2)} = R_{B}D_{B}R_{B}^{-1}, 
\end{equation}   
where $R_{A}$ and $R_{B}$ are the matrices whose columns are the eigenvectors, and 
$D_{A}$ and $D_{B}$ are diagonal matrices of the eigenvalues of $A$ and $B^T_{(2)}$. The matrices $A$ and $B^T_{(2)}$ admit the spectral decomposition \eqref{eq_spectral} because the multiplicities of their eigenvalues are equal to one. This property can be established at least numerically, see \cite{SMNF, SemDiss}.

We multiply the equation \eqref{eq_Iter} from the left by $R_{A}^{-1}$ and from the right by $R_{B}$ and introduce the notation: $V^{n+1} = R_{A}^{-1} U_0^{n+1} R_{B}$. Then, equation \eqref{eq_Iter} can be written in the form:
\begin{eqnarray} 
 \label{eq_Iter1}
	&&V^{n+1} - \Delta D_{A} V^{n+1}   +  \Delta \nu V^{n+1} D_{B}  =   \\
	&&\underbrace{ R_{A}^{-1} \Big(U_0^{n} + \Delta\left[ (U_0^{n} + U_{C_k})\odot \big((U_0^{n} + U_{C_k}) B_{(1)}^T\big)\right] \Big) R_B - \nu \Delta R_{A}^{-1} U_{C_k} R_{B} D_{B} }_{\text{denote as $G^{n+1}$}}. \nonumber
\end{eqnarray}   

We express the elements of the matrix $V^{n+1}$ standing in \eqref{eq_Iter1} as
\begin{equation}
\label{eq_V}
	V_{ij}^{n+1} = G^{n+1}_{ij} \big/ \left(1 - \Delta D_{A,i} + \Delta \nu D_{B,j}\right),
\end{equation} 
where $D_{A,i}$ is the $i$-th eigenvalue of the matrix $A$, $D_{B,j}$ is the $j$-th eigenvalue of $B^T_{(2)}$, $G^{n+1}_{ij}$ is the value of the right-hand side of \eqref{eq_Iter1} at the node $(t_i,x_j)$. We express the solution at the next iteration of the relaxation method in the form  $U_0^{n+1} = R_{A} V^{n+1} R_{B}^{-1}$, increase $n$ by unity and repeat the described computations.

We perform iterations of the relaxation method in the current cell until the solution stabilizes. As a criterion for stopping the iterations, we use the condition
\begin{equation}
\label{residd}
\max|U_0^{n+1}-U_0^n| \leq \varepsilon_S \Delta,
\end{equation}
where $\varepsilon_S$ is the residual of relaxation, which is assumed to be small.
Here, the maximum is taken over all elements of the matrix.
In this algorithm we set $U^{0}_{0}$ to be the zero matrix. If the stopping criterion is not satisfied in the current cell  $C_k$ for $n = \overline{1,199}$, then we reduce the time step $d_k$ by half and rerun the computations in the cell $C_k$. We continue these operations until the relaxation method converges in cell $C_k$ with fewer than 200 iterations.

We can find the numerical solution $u(\tau_{k+1},x)$ using \eqref{eq_int}. Following the ideas of \cite{TeeTref}, we describe the method for localizing the singularity of the solution $z^*(t) = \delta(t) + i \varepsilon(t) \in\mathbb{C}$ and for adapting the spatial grid $x_j$, $j=\overline{1,N}$. To implement this method, we assume that the solution has an isolated singularity lying close enough to the domain of the problem.
To find approximate values of $\delta(\tau_{k+1})$ and $\varepsilon(\tau_{k+1})$, we expand $u(\tau_{k+1},x)$ in a Fourier-Chebyshev series and then apply a Chebyshev--Pade approximation to it:
\begin{equation} 
\label{eq_pade_cheb}
	u(\tau_{k+1},x) = \sum_{k=0}^{\infty} a_k T_k(x)  \approx \dfrac {\sum_{k=0}^m b_k T_k(x)} {1+\sum_{k=1}^r c_k T_k(x) } = \dfrac {p(x)}{q(x)}.
\end{equation} 
By requiring that $u(\tau_{k+1},x) q(x) - p(x) = O(T_{m+r+1}(x))$, we obtain the system of equations: 
\begin{equation*} 
\begin{pmatrix}
1 &           &    &  -d_{0,1}     & \cdots & -d_{0,r}    \\
  & \ddots &    &     \vdots      & \vdots & \vdots        \\
  &            & 1 &  -d_{m,1}    & \cdots &-d_{m,r}    \\
  &            &    &  -d_{m+1,1}& \cdots &-d_{m+1,r}\\
  &            &    &  \vdots         & \vdots & \vdots        \\
  &            &    &  -d_{m+r,1}& \cdots &-d_{m+r,r}\\
\end{pmatrix}
\begin{pmatrix}
b_0    \\
\vdots \\
b_m    \\
c_1     \\
\vdots \\
c_r     \\
\end{pmatrix}
=
\begin{pmatrix}
a_0        \\
\vdots    \\
a_m       \\
a_{m+1}\\
\vdots     \\
a_{m+r}\\
\end{pmatrix}.
\end{equation*}
The values of $d_{i,j}$ can be derived using the condition $\sum_{i=0}^{\infty} d_{i,j} T_i(x) = \sum_{i=0}^{\infty} a_i T_i(x) T_j(x)$ and the formula $T_i(x)T_j(x)=\frac{1}{2}(T_{i+j}(x)-T_{|i-j|}(x))$. To localize the singularity, we need to find $c_i$, $i = \overline{1,r}$. Since the last $r$ equations of the above system do not include $b_i$, $i =\overline{0,m}$, we need to solve only the subsystem for the coefficients $c_i$. If we take $r = 2$, then these coefficients can be expressed explicitly \cite{Bryndin}: first, we need to express
\begin{equation*} 
\begin{matrix}
d_{m+1,1} = a_{m+1}/2 +a_{m+3}/2, & d_{m+1,2} = a_{m}/2 +a_{m+4}/2, \\
d_{m+2,1} = a_{m+2}/2 +a_{m+4}/2, & d_{m+2,2} = a_{m+1}/2 +a_{m+5}/2.
\end{matrix}
\end{equation*}
Then, after solving the subsystem, we need to find the roots of the denominator of the approximation \eqref{eq_pade_cheb}. If the roots are a pair of complex conjugate numbers and the moduli of their imaginary parts are greater than $10^{-14}$, then we take $\delta(\tau_{k+1})$ and $\varepsilon(\tau_{k+1})$ to be the real and imaginary parts of one of the roots. Without loss of generality, we choose the one with positive $\varepsilon(\tau_{k+1})$. If $|\varepsilon(\tau_{k+1})|\leq 10^{-14}$ or if the moduli of $\delta(\tau_{k+1})$ and $\varepsilon(\tau_{k+1})$ are large enough, then we set $\delta(\tau_{k+1})=\delta(\tau_{k})$ and $\varepsilon(\tau_{k+1})=\varepsilon(\tau_{k})$.

To find the coefficients $a_k$ in \eqref{eq_pade_cheb}, we can compute $u(\tau_{k+1},x)$ at the Chebyshev nodes in the segment $[-1,1]$, and then perform a fast Fourier transform. 

If the location of singularity at the previous time step $\tau_{k}$ is known, then, to compute the values $\varepsilon(\tau_{k+1})$, $\delta(\tau_{k+1})$ more accurately, we construct the approximation \eqref{eq_pade_cheb} on the segment $[\max\big(-1,\delta(\tau_k) - 10\varepsilon(\tau_k)\big), \min\big(1,\delta(\tau_k) + 10\varepsilon(\tau_k)\big)]$.

Finally, once $\delta(\tau_{k+1})$ and $\varepsilon(\tau_{k+1})$ are obtained, we can define a new spatial grid for the next cell $C_{k+1}$. For this purpose, we use the mapping $g(x)$ that conformally maps the segment $[-1,1]$ onto itself (see \cite{TeeTref}):
\begin{equation} 
\label{gMap}
\hat x_j = g(x_j) = \delta(\tau_{k+1}) + \varepsilon(\tau_{k+1}) \sinh \Big[ (a^{-} + a^{+})\dfrac{ x_j - 1}{2} + a^{-} \Big],
\end{equation}
where $a^{\pm} = \sinh^{-1}\dfrac{1\pm\delta(\tau_{k+1})}{\varepsilon(\tau_{k+1})}$, $j = \overline{1,N}$.

An important observation is that obtaining a numerical solution to the problems \eqref{eq_ub}--\eqref{eq_ub2} using \eqref{gMap} turns out to be impossible for $\nu<10^{-5}$. In this regard, we develop a modification of \eqref{gMap}. An analysis of the solution to \eqref{eq_ub}, \eqref{eq_ub1} shows that the real coordinate of the singularity can be approximated by zero with high accuracy, see Fig.~\ref{Bfig1}a, therefore we set $\delta(t) = 0$ and instead of $g(x)$ we use the composition $g\circ g$:
\begin{equation} 
\label{gcg}
\hat x_j = g(g(x_j))= \varepsilon(\tau_{k+1}) \sinh \Big[  \frac{\pi}{2} \sinh \Big( x_j S(\tau_{k+1}) \Big)\Big],~j = \overline{1,N},
\end{equation}
where $S(\tau_{k+1})=\sinh^{-1}\Big(\frac{2}{\pi}\sinh^{-1}\frac{1}{\varepsilon(\tau_{k+1})}\Big)$. This composition can take the singularity $z^*(t)$ farther away from the segment $[-1,1]$ in $\mathbb{C}$ than the original mapping $g$. The use of $g\circ g$ in the case $\nu\geq 10^{-5}$  allowed us to significantly improve the stability of the algorithm and to obtain higher accuracy for a fixed number of nodes. As written above, for $\nu<10^{-5}$ the use of $g\circ g$ was the only way to compute the solution.

To improve the stability of the algorithm in the case of the problem \eqref{eq_ub}, \eqref{eq_ub2} with $\nu \leq 10^{-4}$, we use the composition $g\circ g$ presented in the form
\begin{equation}
\label{gMap2}
\hat x_j = g(g(x_j)) =\delta(\tau_{k+1}) + \varepsilon(\tau_{k+1})\Pi\Big[a^{\pm},\delta^*(\tau_{k+1}) + \varepsilon^*(\tau_{k+1}) \Pi(b^{\pm},x_j)\Big],~j = \overline{1,N},
\end{equation}  
where 
$$
\Pi(\alpha^{\pm},x) = \sinh\Big[(\alpha^++\alpha^-)\frac{x-1}{2}+\alpha^-\Big],~~~\delta^*(\tau_{k+1})  = \dfrac{a^+-a^-}{a^++a^-},$$
$$
\varepsilon^*(\tau_{k+1}) = \dfrac{\pi}{a^++a^-},~~~b^{\pm} = \sinh^{-1}\dfrac{1\pm\delta^*(\tau_{k+1})}{\varepsilon^*(\tau_{k+1})}.
$$

\section{Numerical results and discussion}
\label{S4}

In this section, the results of studying the accuracy and time performance of the constructed algorithm are presented in the case of solving initial-boundary value problems \eqref{eq_ub}--\eqref{eq_ub2} for the Burgers equation. In sect.~\ref{secc1} we describe a method for computing the exact solution \eqref{eq_ub}, \eqref{eq_ub1} by formula \eqref{eq_ub_exact} and study its convergence. In what follows, we use this method to verify the numerical results. In sect.~\ref{secc2} it is shown that even for $\nu=0.01$ the construction of adaptive methods with a justified choice of the adaptation is extremely important for computations. Then, implementation details of the proposed method for problem \eqref{eq_ub}, \eqref{eq_ub1} along with an analysis of the obtained solutions are presented. In sect.~\ref{secc3} these solutions are compared with the ones obtained by other authors (we perform a comparison of the values of the derivative $\partial u/\partial x$ at $x=0$). To evaluate the quality of the results, we compute the exact value of this derivative using formula \eqref{eq_ub_exact}. It should be emphasized that the solutions of the considered problem for $\nu<10^{-3}$ were not found in literature. However, we were able to compute solutions for extremely small $\nu$ down to $\nu=10^{-8}$. In sect.~\ref{secc4} the proposed algorithm based on the STSM is compared with an algorithm that uses high-order Runge--Kutta scheme, RK5(4) 7S, for time stepping, which is described in \cite{Dormand-Prince} and implemented in \cite{Bryndin}. Furthermore, we characterize the cases in which the STSM is more effective and those in which RK5(4) 7S demonstrates higher efficiency. Finally, in sect.~\ref{secc5}, we study the problem \eqref{eq_ub}, \eqref{eq_ub2}, which has the exact solution \eqref{exSol}, and compare our numerical solutions with the results from \cite{Dogan}. It is demonstrated that the proposed method enables us to achieve significantly higher accuracy at a low computational cost. Furthermore, we compute the approximate trajectories of the singular points \eqref{sing_exSol} of solution \eqref{exSol} lying in a small neighbourhood of the domain of problem. We compare these trajectories (the values of $\delta(t)$ and $\varepsilon(t)$ for $t\in[T_0,T]$) with the exact ones. It is shown that for values of $\nu$ down to $10^{-6}$ they are computed with high accuracy.

\subsection{Computation of the exact solution}
\label{secc1}

To compute the exact solution of the problem \eqref{eq_ub}, we calculate the integrals appearing in formula \eqref{eq_ub_exact} at time $t=T$ using the double exponential quadrature rule described in \cite{Takahasi}:
\begin{equation}
\label{eq_Ih}
I_h = \dfrac{\pi}{2}h \sum_{n=-H}^{H}w\Big(t,x,\sinh{(\dfrac{\pi}{2}\sinh{n h})}\Big)(\cosh{nh})\cosh{\Big(\dfrac{\pi}{2}\sinh{nh}\Big)},
\end{equation}
where $h=  \sinh^{-1} \left(\dfrac{2}{\pi} \sinh^{-1}(10^{19})\right) \Big/ H$, $2H+1$ is the number of nodes of the quadrature rule, and $w(t,x,y)$ is the integrand. The number $10^{19}$ in the expression for the step $h$ serves as a numerical approximation of infinity. When computing \eqref{eq_ub_exact} with the initial condition $\varphi(x)=-\sin \pi x$ and the viscosity $\nu < 10^{-3}$, we use a correction factor with function $\psi(t,x)$ calculated as described in sect.~\ref{S2}. 

\begin{figure}[h!]
\centering 
\includegraphics[width=0.40\linewidth]{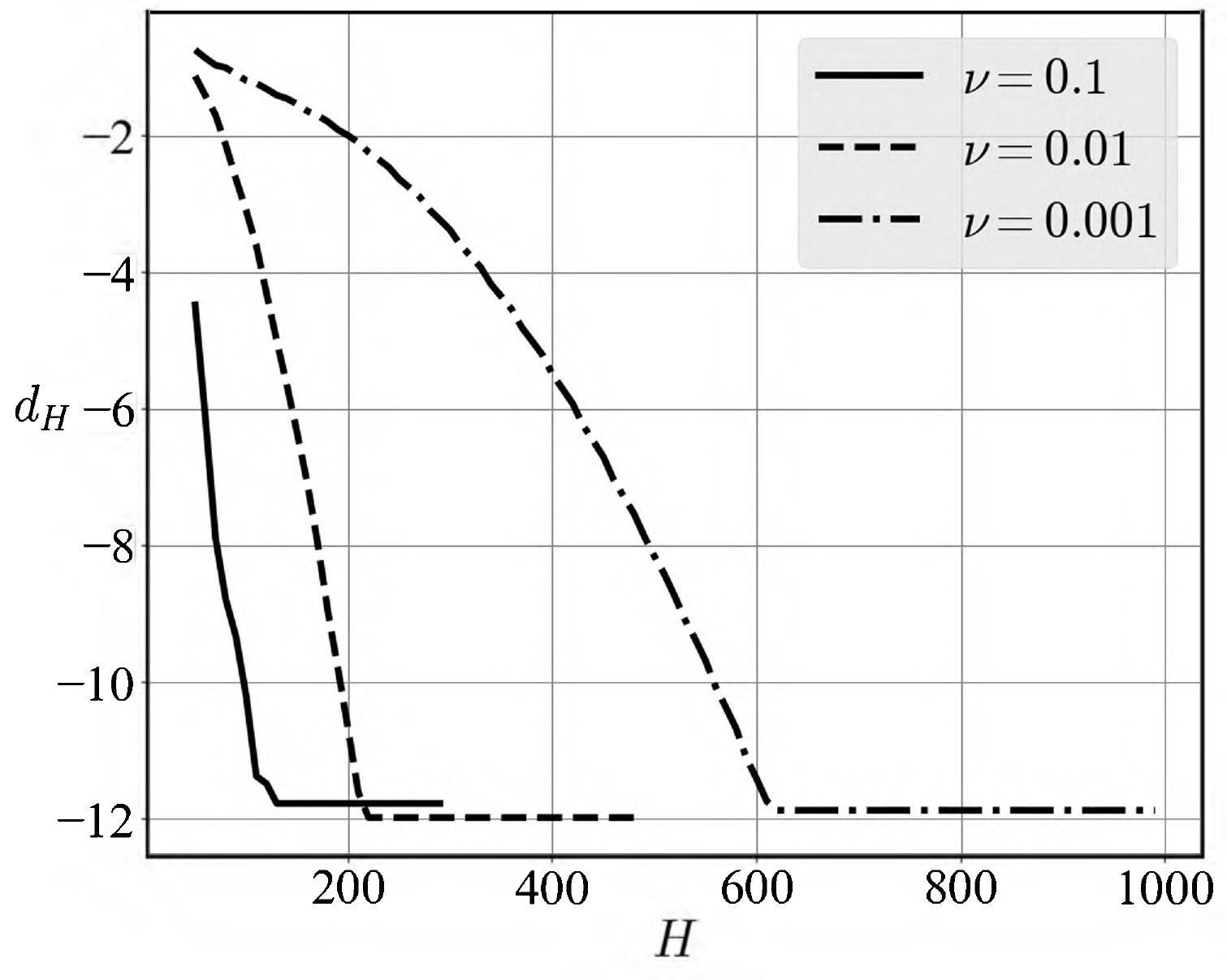}
\caption{Plot of $d_{H}$ as a function of $H$ at $T =0.5$.}
\label{figErr1}
\end{figure}

Fig.~\ref{figErr1} shows the convergence of the approximations of the exact solution \eqref{eq_ub_exact}, where the integrals entering \eqref{eq_ub_exact} are computed by formula \eqref{eq_Ih}. Here, the case of $\varphi(x)=-\sin\pi x$ and of various values of $\nu$ is considered. The plot shows a quantity $d_{H} = \lg\left(\dfrac{\|u_{ex,H}(T,x)-u_{ex,2H}(T,x)\|_{\infty}} {\|u_{ex,2H}(T,x)\|_{\infty}}\right)$ as a function of $H$ at $T =0.5$. Here, $\|\cdot\|_{\infty}$ is the supremum norm, which is computed as described below; $u_{ex,H}(T,x)$ are the values of the exact solution computed via \eqref{eq_ub_exact} employing a double exponential quadrature rule with $2H+1$ nodes. We can conclude that, for sufficiently large $H$, the solution $u_{ex,H}$ is highly accurate. Note that for smaller values of $\nu$ we did not construct such plots, but we always verified that doubling $H$ changes the relative error $err$ (defined below) by no more than 1\%.

\subsection{Result for the problem \eqref{eq_ub}, \eqref{eq_ub1}}
\label{secc2}

To study the accuracy of the computations, we define a uniform grid over $x$ with $N_p$ nodes: $\tilde X=\{x_j: x_j=-1+2(j-1)/(N_p-1), j=\overline{1,N_p}\}$ and introduce the notation for the relative error of the numerical solution to problem \eqref{eq_ub}:
\begin{equation}
\label{Rerr}
err =  \dfrac{\|u(T,x) - u_{ex}(T,x)\|_\infty}{\|u_{ex}(T,x)\|_\infty},
\end{equation}
where $\| \cdot \|_\infty = \max\limits_{x_j\in\tilde X} |\cdot|$, $u(T,x)$ is the value of the approximate solution computed at time $t=T$ at the point $x$, $u_{ex}(T,x)$ is the value of the exact solution given either by an explicit formula or computed as described in the previous section. If no value of $N_{p}$ is specified, it is assumed that $N_p = 10000$. 

For a more detailed study of the accuracy of the numerical solution to problem \eqref{eq_ub}, we also define the local relative error
\begin{equation*}
err_{loc} =  \dfrac{\|u(T,x) - u_{ex}(T,x)\|_{\infty, loc}}{\|u_{ex}(T,x)\|_{\infty, loc}},
\end{equation*}
where $\| \cdot \|_{\infty, loc} = \max\limits_{x_j\in\bar X} |\cdot|$ and $\bar X=\left\{x_j: x_j = -\mu+2\mu\dfrac{(j-1)}{(N_{p, loc}-1)}, j=\overline{1,N_{p, loc}}\right\}$, $\mu = 100\nu$, $x_j \in [-100\nu, 100\nu]$. 

Below, we outline the basic parameters of the numerical method and describe their values used in the computations:
\begin{itemize}
\item[1.] $N$ is the number of nodes of the spatial grid, $N\in[100, 350]$.
\item[2.] $K$ is the number of cells, into which the problem domain is decomposed; it can be fixed or vary during the computation to adapt the algorithm, $K\in[10, 2\times 10^6]$.
\item[3.] $M$ is the number of Chebyshev nodes of the temporal grid in each cell; it must be sufficiently large to provide high accuracy of the approximation \eqref{eq_intetpolate}, but not too large so as not to increase the computational cost and to preserve the numerical stability of the algorithm. In numerical tests we established that $M = 10$ allows one to obtain the minimal $err$.
\item[4.] The size of the step of the relaxation method, $\Delta  = 1$, provides fast convergence in most computational regimes.
\item[5.] $\varepsilon_S$ is the relaxation residual, see \eqref{residd}. Unless otherwise specified, $\varepsilon_S=10^{-12}$.
\item[6.] $H$ is the parameter of the quadrature rule \eqref{eq_Ih} used to compute the exact solution \eqref{eq_ub_exact}, $H\in[300,10^5]$.
\item[7.] $N_{p}$ and $N_{p, loc}$ are the numbers of nodes of a uniform grid in $x$ on which the relative error and the local relative error are computed, respectively, $N_p\in[10^4 ,10^5]$, $N_{p,loc}\in[5\cdot10^3 ,2\cdot10^5]$. Unless otherwise specified, $N_p=10^4$.  
\end{itemize}

\begin{figure}[h!]
\centering 
\includegraphics[width=0.4\linewidth]{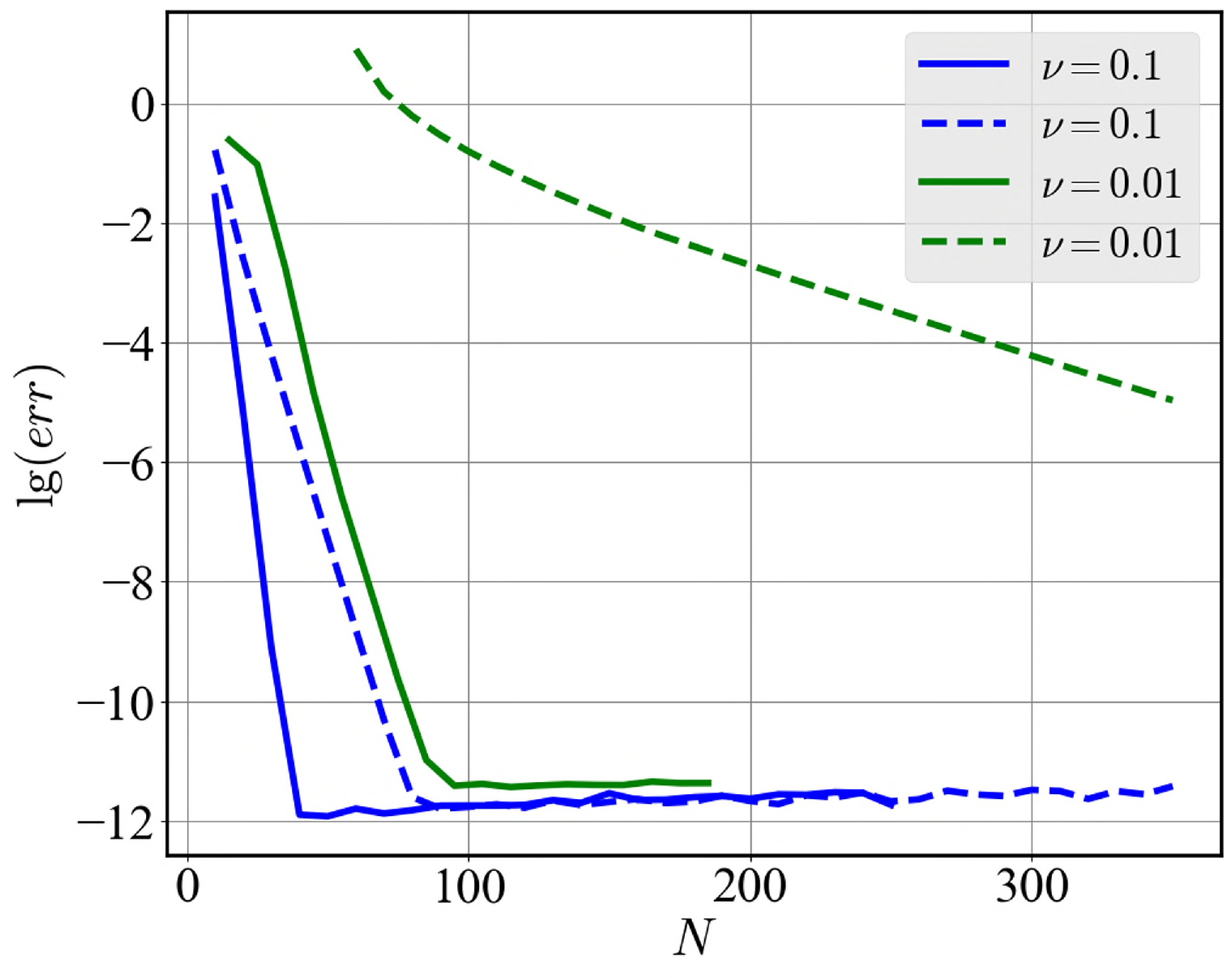} 
\caption{ Plots of $\ln(err)$ as a function of $N$ for $T =0.5$ and $\nu=0.1, 0.01$ in the case of the adaptive algorithm (solid lines) and the algorithm without adaptation (dashed lines). Parameters of the algorithm in the case $\nu = 0.1$ are $M = 10 $, $K=10 $, $H = 300$; in the case $\nu = 0.01$ -- $M = 10 $, $K=20 $, $H = 300$. }
\label{figErr2}
\end{figure}

The application of RA with respect to the spatial variable and of the adaptation of the grid to the location of the singularity allow one to obtain a highly accurate numerical solution for small values of $\nu$. In Fig.~\ref{figErr2}, the plots of $\lg(err)$ versus $N$ for $T =0.5$ and $\nu=0.1, 0.01$ are presented for the solutions obtained using the algorithm with adaptation of the spatial grid to the singular point (solid lines) and without adaptation (dashed lines). In these experiments, the mapping $g$ \eqref{gMap} was applied to adapt the grid to the singularity. The solid lines on the presented plot show a fast exponential decrease of the error for sufficiently small $N$. To achieve a comparable accuracy using the algorithm without adaptation, a considerably larger number of nodes $N$ is required, which characterizes the efficiency of the adaptive method. Its efficiency is also confirmed by the fact that, if the value of viscosity is less than $10^{-2}$, the algorithm without adaptation of the grid does not converge at all. Adaptation with a single application of the map $g$ \eqref{gMap} allowed us to obtain the numerical solutions for $\nu\gtrsim 10^{-5}$, while the use of the composition $g\circ g$ described in \eqref{gcg} allowed us to decrease the values of $\nu$ down to $10^{-8}$. 

To test the adaptive algorithm, we conduct a number of numerical experiments, the results of which are summarized in Tables~\ref{SINtabl0.5} and \ref{SINtabl1}. In all of these simulations, the composition $g\circ g$ was used, and the values of the parameter $N_p$ were chosen so that doubling $N_p$ changes the values of $err$ and $err_{loc}$ by no more than 1\%. 
For some of the numerical solutions, the values of the error attain their maxima only in a neighbourhood of $x=0$. In such cases, the accuracy is characterized by $err_{loc}$. In Tables~\ref{SINtabl0.5} and \ref{SINtabl1}, the values of $\max(err,err_{loc})$ are given. Table~\ref{SINtabl0.5} includes the results of eight experiments, and Table~\ref{SINtabl1} includes seven experiments; the experiments (tests) are numbered. In the first two tests described in each table, the number of cells $K$ is fixed.
When performing tests numbered 3 and higher, the initial size of the computational cell was set to $d_k = 0.005$. In test No.~3 at $t = 0.4375$, the relaxation method performed 200 iterations. As a result, in accordance with the recommendations given after \eqref{residd}, the size of the cell was halved and the solution was recomputed. A similar situation occurred in all higher-numbered tests at various $t$.

\begin{table*}[htb!]
\caption{Testing the STSM by solving problem~\eqref{eq_ub},~\eqref{eq_ub1} on the segment $t\in[0,T]$, $T=0.5$ for various $\nu$. Here $time$ is the running time (s are seconds, m are minutes, h are hours).}
\begin{center}
\begin{tabular}{l|l|l|l|l|l|l|l} 
	       & $\nu$   & $N$   & $K$    & $H$     & $\max(err, err_{loc})$ & $time$   & $\varepsilon(T)$ \\
 \hline 1 &  1.e-2  & 100   & 10		& 3e+02   & 1.82e-12  & 14s	  &  3.37e-02	\\   
 \hline 2 &  1.e-3  & 150   & 20		& 1e+03   & 4.57e-12  & 1m 30s	  &  3.15e-03 \\    
 \hline 3 &  1.e-4  & 150   & 125	& 5e+03   & 1.17e-10  & 9m 35s	  &  3.14e-04 \\     
 \hline 4 &  1.e-5  & 300   & 693	& 2e+04   & 4.89e-09  & 6h 29m   &  3.14e-05 \\  
 \hline 5 &  1.e-5  & 250   & 638	& 2e+04   & 5.73e-07  & 1h 20m   &  3.14e-05 \\ 
 \hline 6 &  1.e-6  & 200   & 5936	& 2.5e+04 & 3.53e-07  & 4h 50m   &  3.14e-06 \\ 
 \hline 7 &  1.e-7  & 250   & 48988 	& 5e+04   & 9.18e-05  & 27h 	  &  3.14e-07 \\ 
 \hline 8 &  1.e-8  & 250   & 450441 & 1e+05   & 6.83e-03  & 72h       &  3.44e-08 
\end{tabular}
\end{center}
\label{SINtabl0.5}
\end{table*} 
\begin{table}[htb ]
\caption{Testing the STSM by solving problem~\eqref{eq_ub},~\eqref{eq_ub1} on the segment $t\in[0,T]$, $T=1$ for various $\nu$. Here $time$ is the running time (s are seconds, m are minutes, h are hours).}
\begin{center}
\begin{tabular}{l|l|l|l|l|l|l|l} 
	       & $\nu$   & $N$   & $K$      & $H$     & $\max(err, err_{loc})$ & $time$  & $\varepsilon(T)$ \\ 
 \hline 1 &  1.e-2  & 100   & 10		  & 3e+02   & 1.77e-12  & 16s	    &  4.31e-02 \\      
 \hline 2 &  1.e-3  & 150   & 40	 	  & 1e+03   & 8.66e-12  & 3m 02s	&  4.27e-03 \\      
 \hline 3 &  1.e-4  & 150   & 325	  & 5e+03   & 2.57e-10  & 24m 54s	&  4.26e-04 \\  
 \hline 4 &  1.e-5  & 250   & 2658	  & 2e+04   & 5.43e-07  & 3h 44m 	&  4.26e-05 \\  
 \hline 5 &  1.e-6  & 200   & 25936	  & 2.5e+04 & 3.73e-06  & 21h    	&  4.26e-06 \\ 
 \hline 6 &  1.e-7  & 250   & 208988   & 5e+04   & 2.93e-04  & 48h 	    &  4.26e-07 \\ 
 \hline 7 &  1.e-8  & 250   & 1730442  & 1e+05   & 2.02e-02  & 336h	    &  4.33e-08 
\end{tabular}
\end{center}
\label{SINtabl1}
\end{table}

The constructed algorithm makes it possible to achieve a good balance between accuracy and time performance that is important for solving applied problems. This balance can be achieved by changing the parameters $N$ and $K$: their reduction leads, on the one hand, to a loss of accuracy, and on the other hand, to a decrease in the running time of the algorithm. Moreover, to accelerate the algorithm --- admittedly at the cost of reduced accuracy --- we can use the following techniques (see the results of tests No.~5--8 in Table~\ref{SINtabl0.5} and No.~4--7 in Table~\ref{SINtabl1}):

\begin{itemize} 
    \item[1)] If $ 10^{-6}\leq \nu < 10^{-4} $, then reduce $d_k$ to $0.001$ after $t = 0.3$; perform localization of the singularity, refinement of the grid and spectral decomposition of $B^T_{(2)}$ no more than once over a time segment of length $0.001$;
    \item[2)] If $\nu < 10^{-6}$, then reduce $d_k$ to $0.0001$ after $t = 0.3$; perform localization of the singularity, refinement of the grid and spectral decomposition of $B^T_{(2)}$ no more than once over a time segment of length $0.0001$.
\end{itemize}

\begin{figure}[h!]
\centering
\includegraphics[width=0.4\linewidth]{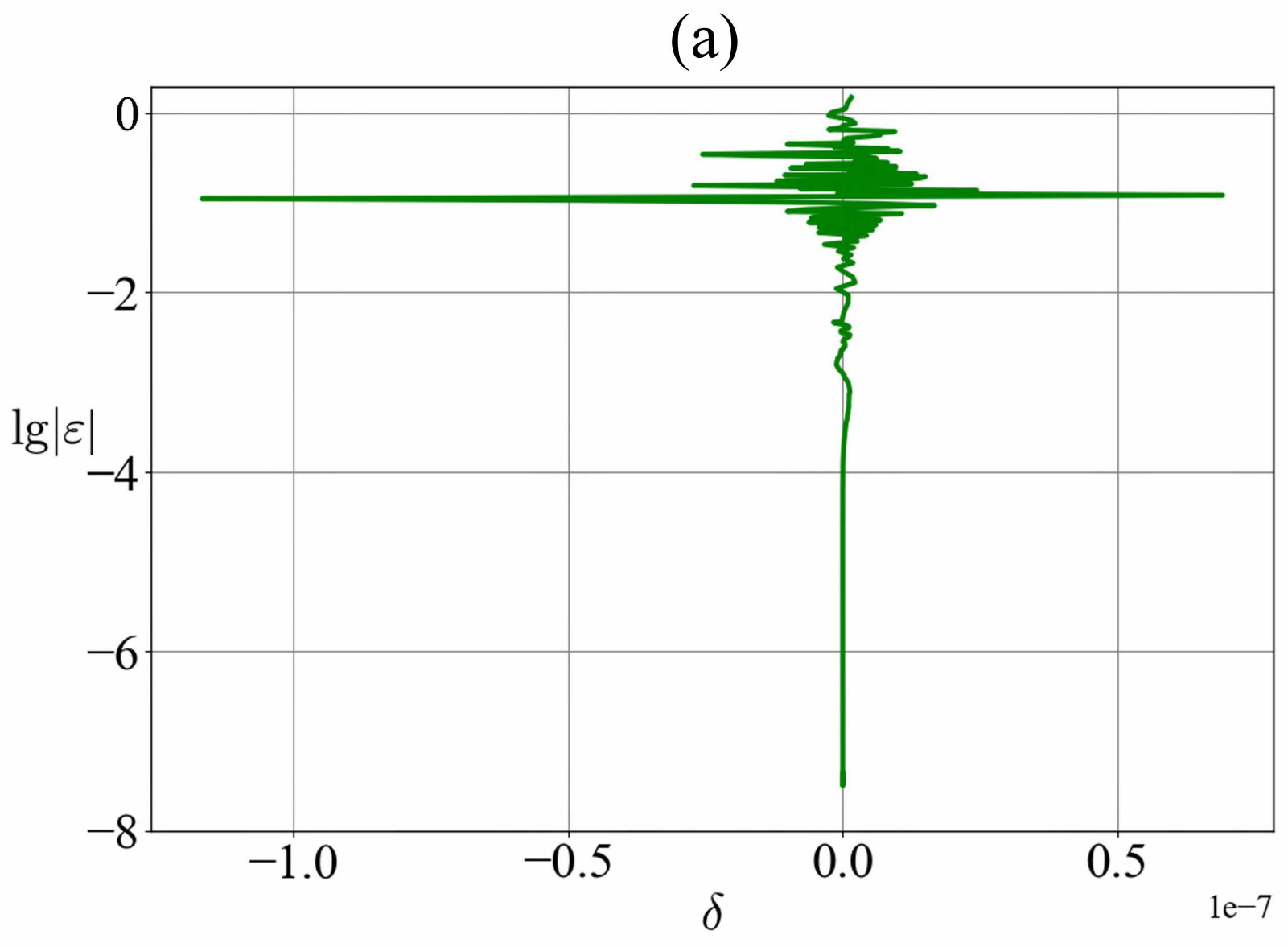}
\includegraphics[width=0.4\linewidth]{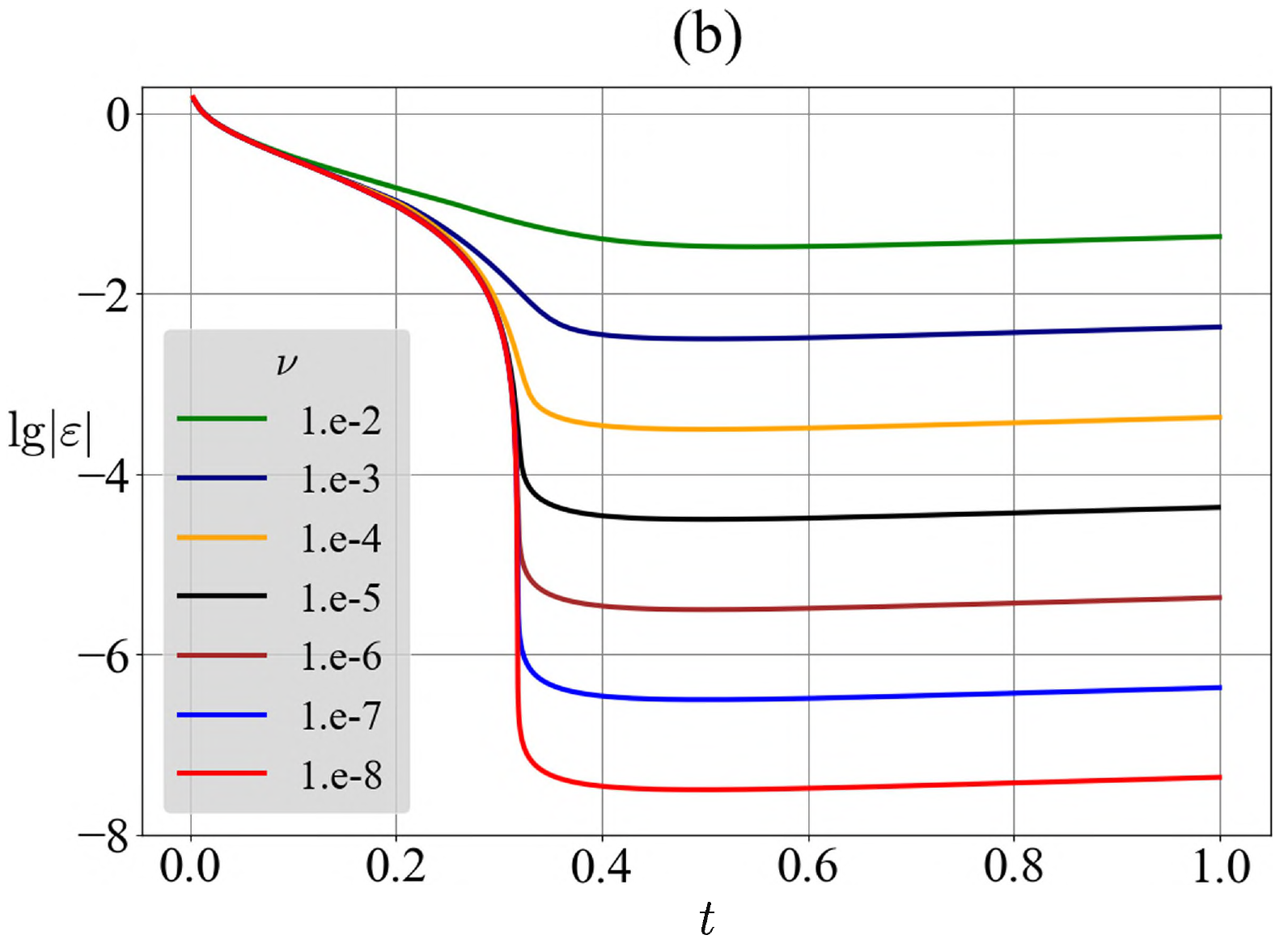}
\caption{ Trajectory of the singular point in problem \eqref{eq_ub}, \eqref{eq_ub1}, $t \in [0,1]$ for $\nu = 10^{-8}$, $N = 250$, $M=10$ and $K = 1730442$ (a); dependencies of the imaginary parts of the coordinates of the singular points on $t$ for various values of $\nu$ (b).}
\label{Bfig1}
\end{figure}

At time $t \approx 0.4$, the imaginary part of the coordinate of the singular point reaches a minimum (Fig.~\ref{Bfig1}b), which leads to the emergence of large gradients of the solution with respect to the variables $x$ and $t$. To approximate functions with such gradients in $t$, one has to either increase the number of nodes $M$ or reduce the grid step $d_k$. Our experiments show that the second option (reducing $d_k$) yields considerably better results.

We should also note that for a large value of $N$, the spectral decomposition of $B^T_{(2)}$ in each cell takes a considerable time. Therefore, to accelerate the algorithm, it is proposed to skip this procedure in a small neighbourhood of the time moment at which the singular point has already been located and the grid adapted. This suggestion certainly leads to a loss of accuracy but makes it possible to considerably accelerate the algorithm. This can be observed in tests 4--5 from Table~\ref{SINtabl0.5}. In test No.~4, neither of the two acceleration techniques described above was used, which yielded a more accurate solution, but the running time was about 5 times greater than in test No.~5, where both recommendations were applied.

The application of the composition $g\circ g$ for $\nu \geq 10^{-5}$ gave a much smaller error than the use of the original map $g$, see \eqref{gMap}. For $\nu < 10^{-5}$, the application of this composition turned out to be the only way to make the algorithm converge.

From the results given in Tables~\ref{SINtabl0.5} and \ref{SINtabl1}, we can derive several patterns. First, for sufficiently large $t$, $\varepsilon(t)\approx C_{t}\nu$, where $C_{t}$ has a slow growth as $t$ increases. Second, in the case of adapting the step size $d_k$, the number of cells $K$ is approximately inversely proportional to $\nu$. Third, computations on the segment $t\in[0.5,1]$ usually require significantly more time than computations on the segment $t\in[0,0.5]$. This is naturally explained by the fact that for all $t\in[0.5,1]$ the singular point lies close to the domain of problem, see Fig.~\ref{Bfig1}. Moreover, if both methods of accelerating the algorithm described above are used, the error of the solution at $t=0.5$ is approximately an order of magnitude lower than at $t=1$.

In Fig.~\ref{Bfig1} the trajectory of the singular point that is closest to the domain of the problem is plotted. This trajectory is computed applying the numerical algorithm described above, which uses the composition of maps $g\circ g$. The figure shows that the imaginary coordinate of the singular point initially rapidly approaches its minimum value until the time $t^* \approx 0.4$ (this value depends on $\nu$) and then bounces off from this value, but the order of magnitude of $\varepsilon(t)$ for $t>t^*$ is preserved. For small $t$, fluctuations of the real part $\delta(t)$ of the coordinate of the singular point are observed. This effect probably occurs due to the numerical error of the approximate method. It does not influence the computational process much because of the negligibly small amplitude ($\sim 10^{-7}$) of these fluctuations, see Fig.~\ref{Bfig1}a. We also note that the dynamics of $\varepsilon(t)$, shown in Fig.~\ref{Bfig1}b, qualitatively corresponds to the data given in \cite{SulemSulemFrish} for $\nu=0.05$. This dynamics is inherent to the solution of the BE with small $\nu$ and the initial condition representing an entire function: at the initial moment, the singular points of the solution come from infinity, rapidly approach the segment of the problem, and then begin to move away from it (more information about this effect is given in \cite{BE_CP}).

\subsection{Comparison of the obtained solutions to problem \eqref{eq_ub}, \eqref{eq_ub1} with the results other authors}
\label{secc3}

The problem \eqref{eq_ub}, \eqref{eq_ub1} was studied in a number of papers, see~\cite{SulemSulemFrish,Basdevan,OVasilyev} and the references therein. However, because of difficulties in computing the exact solution (see sect.~\ref{S2}), the numerical result could not be compared with it for small values of $\nu$. Nevertheless, values of the derivative $\partial u/\partial x$ at $x=0$ obtained for moderate values of $\nu$ and at different moments of time are available in several publications. Now, we can see that in those works they were computed with a rather large error.

In Table~\ref{Gradient}, the values of the mentioned derivative are presented at time $T=0.5$ for various values of $\nu$. The exact values of the derivative are computed by differentiating the solution \eqref{eq_ub_exact}, applying a correction factor with the function $\psi(t,x)$ and using the double exponential quadrature rule with $H = 10^5$ and $2H$. The doubling of $H$ enables us to control the error. The exact value of the derivative for given $H$ and $T$ is denoted by $$D_{ex}^{H,T}=\frac{\partial u_{ex}}{\partial x}\big|_{x=0}.$$
In the second column of Table~\ref{Gradient}, the last four digits of the value of $D_{ex}^{2H,T}$ are given in parentheses, from the 13th to the 16th digit. This shows that the values of the derivative are computed with an accuracy of at least 13 digits.

In the third and fourth columns of Table~\ref{Gradient}, the respective approximate values of the derivatives are presented. They are obtained by taking the most accurate results from Table~\ref{SINtabl0.5} and differentiating them numerically. Two approaches of the differentiation are implemented: in the first, the barycentric formula without modification is used, see \eqref{bar0}, \eqref{TeeTreff1}. In this case, the value of the derivative at the point $x=0$ is denoted by $D_{WM}$. In the second approach, the modified formulas \eqref{eq_baricentr} and \eqref{eq_diffX1} are used. In this case, the value of the derivative is denoted by $D_{M}$. 

Table~\ref{Gradient} shows that for moderate $\nu$ ($\nu>10^{-5}$), the value of $D_{M}$ is less accurate than that of $D_{WM}$. This is due to the fact that the use of modification \eqref{eq_baricentr} slightly spoils the numerical stability of the approximation method (we postulate without proof that the Lebesgue constant of approximation \eqref{eq_baricentr} grows faster with increasing $N$ than that of \eqref{bar0}). However, this modification is necessary to resolve boundary conditions in the framework of the proposed algorithm, see remark~\ref{rem3}. Therefore, to find a numerical solution of the BE, it is necessary to use the formulas with modification, but to compute the derivative of the solution it is better to use the original barycentric formulas.

For small $\nu$ ($\nu\leq 10^{-5}$), the dominant contribution to the errors of $D_M$ and $D_{WM}$ comes from the error of the numerical solution itself, not from the way of computing the derivative. For such values of $\nu$, the numerical solution has a rather large error and the difference between $D_M$ and $D_{WM}$ becomes smaller than it. Therefore, the values of $D_M$ and $D_{WM}$ practically coincide, and the choice of the method of numerical differentiation no longer affects the accuracy of the result.

\begin{table}[h]
\caption{Values of the derivative of the solution to \eqref{eq_ub}, \eqref{eq_ub1} with respect to $x$ at the point $x=0$ for $T=0.5$ and various $\nu$; the exact values of the derivative are computed for $H = 10^5$ and $2H$.}
\begin{center}
\begin{tabular}{l|l|l|l} 
         $\nu$ & $-D_{ex}^{H,T}$ (digits 13--16 of $-D_{ex}^{2H,T}$) & $-D_{WM}$ & $-D_M$\\ 
 \hline  1.e-2 & 43.88646098024938(4862)                &  43.8864609801           &  43.9       \\
 \hline  1.e-3 & 494.9830870739711(9639)            &  494.9830870732          &  494.982    \\
 \hline  1.e-4 & 4995.057436483183(3162)                 &  4995.057435             &  4995.0573  \\
 \hline  1.e-5 & 49995.06442574773(4775)           &  49995.0643              &  49995.0643 \\
 \hline  1.e-6 & 499995.0651206324(6331)      &  499994.7                &  499994.7   \\
 \hline  1.e-7 & 4999995.065190129(0075)        &  5000686.2               &  5000686.2  \\
 \hline  1.e-8 & 49999995.06519861(963*)         &  49457792.1              &  49457792.1 \\
 \end{tabular}
 \end{center}
\label{Gradient}
\end{table} 

Table~\ref{OVasilyevvsUST} shows a comparison of the values of the derivative at $x=0$ with the results from \cite{OVasilyev}. In this work, the algorithm based on the "Mexican hat wavelet" is presented, and a comparison with the results from the paper \cite{Basdevan} is done, which in turn discusses several approaches to the numerical solution of problem \eqref{eq_ub}, \eqref{eq_ub1}. We state that the algorithm described above gives a significantly higher accuracy in the considered problem than all the mentioned methods.
It is worth noting as well that we use a more precise approach for estimating the error of the solutions based on computing $err$ and $err_{loc}$. 
Indeed, in most papers the deviations between the exact and approximate solutions are computed at the nodes of the grid of the numerical method (an auxiliary fine grid is not introduced); moreover, neighbourhoods of the singular point, where the numerical solution usually has high-frequency local oscillations, are sometimes excluded when computing the error, see~\cite{OVasilyev}. We define the relative error of computing the derivative of the solution at $x=0$ by the formula
$$err_x^T = \dfrac{\Big|D_{ex}^{10^4,T} -  D_{WM} \Big|}{\Big|D_{ex}^{10^4,T}\Big|}.$$
In Table~\ref{OVasilyevvsUST}, the most accurate values of the derivative from \cite{OVasilyev} and the results of our computations obtained using the STSM are presented for the same values of $t$. The values of the relative error are also given. Obviously, the proposed algorithm has significant advantage.

\begin{table}[h]
\caption{Comparison of the results obtained using the STSM with the results from \cite{OVasilyev} for $\nu = 10^{-2}/\pi$, $N=150$, $M=10$, $K=20$.}
\begin{center}
\begin{tabular}{l|l|l|l|l|l} 
        $\pi T$ & $-D_{ex}^{10^4,T}$ & $-\frac{\partial u}{\partial x} |_{x=0}$  \cite{OVasilyev} & $err_x^T$ \cite{OVasilyev} &   $-D_{WM}$, STSM & $err_x^T$, STSM\\
 \hline 1.6030  & 152.00508883277 & 149.59 & 1.59e-02 & 152.0050888329 & 8.49e-13 \\
 \hline 1.6035  & 152.00515616723 & 151.48 & 3.45e-03 & 152.0051561677 & 3.07e-12 \\
 \end{tabular}
 \end{center}
\label{OVasilyevvsUST}
\end{table}  

\subsection{Comparison of the STSM with the Runge--Kutta-based method for problem \eqref{eq_ub}, \eqref{eq_ub1}}
\label{secc4}

In Tables~\ref{RKvsUST1} and \ref{RKvsUST}, we present a comparison of the results obtained with the Spectral Time and Space Method (STSM) described above and with an algorithm employing the RK5(4) 7S scheme for time discretization. This scheme converges with the 5th order and it has a wider region of stability compared with other Runge--Kutta schemes having 4 parameters. In both algorithms, a rational barycentric interpolation (RBI) is applied to approximate the solution with respect to the variable $x$ and the interpolation nodes are adapted using the map \eqref{gMap}. However, in the algorithm based on RK5(4) 7S, the original formula \eqref{bar0} is applied, in which the nodes and weights of the interpolation are taken as $\check x_i = \cos(i\pi/N)$, $i=\overline{0,N}$ and $ \check \omega_i =
\begin{cases}
     (-1)^i / 2, ~i=0,N;\\
     (-1)^i, ~ i = \overline{1,N-1}, 
\end{cases}$, while in the STSM its modification \eqref{eq_baricentr} is applied. The idea of combining the Runge--Kutta method and the RBI was proposed in~\cite{TeeTref, Bryndin}. A detailed description of RK5(4) 7S and of its applications is given in \cite{Dormand-Prince, Bryndin}.

\begin{table}[h]
\caption{Comparison of the algorithm based on the STSM with an algorithm, which uses the RK5(4) 7S, for $T = 0.5$ and fixed $err^* \approx 10^{-10}$.}
\begin{center}
\begin{tabular}{l|l|l|l|l|l|l} 
        $\nu$ & $err_{RK}$ & $err_{STSM}$ & $time_{RK}$ & $time_{STSM}$ & RK: $N,\tau$   & STSM: $N, M, d_k$ \\ 
 \hline 1.e-2 & 2.61e-09   & 1.82e-12    & 20h         &  14s         & 150, 1.e-8         &  100, 10, 0.05       \\
 \hline 1.e-3 & 2.94e-10   & 4.43e-12    & 21h         &  1m 30s      & 150, 1.e-8         &  150, 10, 0.025      \\
 \hline 1.e-4 & 1.86e-10   & 2.02e-10    & 30h         &  16m 18s     & 200, 1.e-8         &  150, 10, 0.0025     \\ 
 \hline 1.e-5 & 8.50e-10   & 2.68e-09    & 30h         &  19h         & 200, 1.e-8         &  300, 10, 0.00025    \\ 
 \end{tabular}
 \end{center}
\label{RKvsUST1}
\end{table} 

\begin{table}[h]
\caption{Comparison of the algorithm based on the STSM with an algorithm, which uses the RK5(4) 7S, for $T = 0.5$ and fixed $err^* \approx 10^{-6}$.}
\begin{center}
\begin{tabular}{l|l|l|l|l|l|l} 
        $\nu$ & $err_{RK}$ & $err_{STSM}$ & $time_{RK}$ & $time_{STSM}$ & RK: $N,\tau$   & STSM: $N, M, d_k$, $\varepsilon_S$ \\ 
 \hline 1.e-2 & 2.59e-06   &  2.88e-06   &  30s        &  13s         & 100, 1.e-5         & 100, 10, 0.05, $10^{-5}$    \\ 
 \hline 1.e-3 & 2.16e-06   &  1.34e-06   &  23s        &  52s         & 70, 1.e-5          & 100, 10, 0,0125, $10^{-5}$  \\ 
 \hline 1.e-4 & 1.72e-06   &  3.86e-06   &  26s        &  9m          & 90, 1.e-5          & 130, 10, 0.0025, $10^{-12}$  \\  
 \hline 1.e-5 & 7.52e-06   &  3.82e-06   &  5m 50s     &  4h          & 100, 1.e-6         & 200, 10, 0.00025, $10^{-6}$ 
 \end{tabular}
 \end{center}
\label{RKvsUST}
\end{table} 

In what follows, we denote by $err_{RK}$ and $err_{STSM}$ the errors \eqref{Rerr} of the numerical solutions obtained by the RK5(4) 7S and the STSM, respectively. Also, let $time_{RK}$ and $time_{STSM}$ denote the running times of these algorithms measured in seconds (s), minutes (m) or hours (h); let $\tau_{RK}$ denote the time step of the RK5(4) 7S scheme.

For a fair comparison of the results obtained by the STSM and RK5(4) 7S, since the time step of the RK5(4) 7S is constant, the cell size $d_k$, $k=\overline{1,K}$, of the STSM was set constant, as well. When obtaining the results presented below, for each algorithm we set a goal of computing the solution with an error less than $err^*$ spending the minimum possible time and we decreased the value of $\nu$ down to $\nu = 10^{-5}$. The last two columns of Tables~\ref{RKvsUST1} and \ref{RKvsUST} contain the values of the parameters of the STSM and the Runge--Kutta-based method (RK) that enable us to achieve this goal.

In the first part of the experiments (see Table~\ref{RKvsUST1}) we set $err^*\approx10^{-10}$. From Table~\ref{RKvsUST1} we see that the algorithm based on the STSM and the relaxation method makes it possible to achieve this accuracy significantly faster than the scheme that uses RK5(4) 7S. However, from the results presented in Table~\ref{RKvsUST}, in obtaining which the accuracy requirements were reduced to $err^*\approx10^{-6}$, we see that the RK5(4) 7S scheme is faster. We also note that in the experiments, whose results are presented in Table~\ref{RKvsUST}, the values $\varepsilon_S$ of the relaxation residual were increased in order to accelerate the STSM. For obtaining the results presented in Table~\ref{RKvsUST1}, we use $\varepsilon_S=10^{-12}$. 

\subsection{Numerical analysis of problem \eqref{eq_ub}, \eqref{eq_ub2} and of its singularities}
\label{secc5}

\begin{figure}[h!]
\centering 
\includegraphics[width=0.55\linewidth]{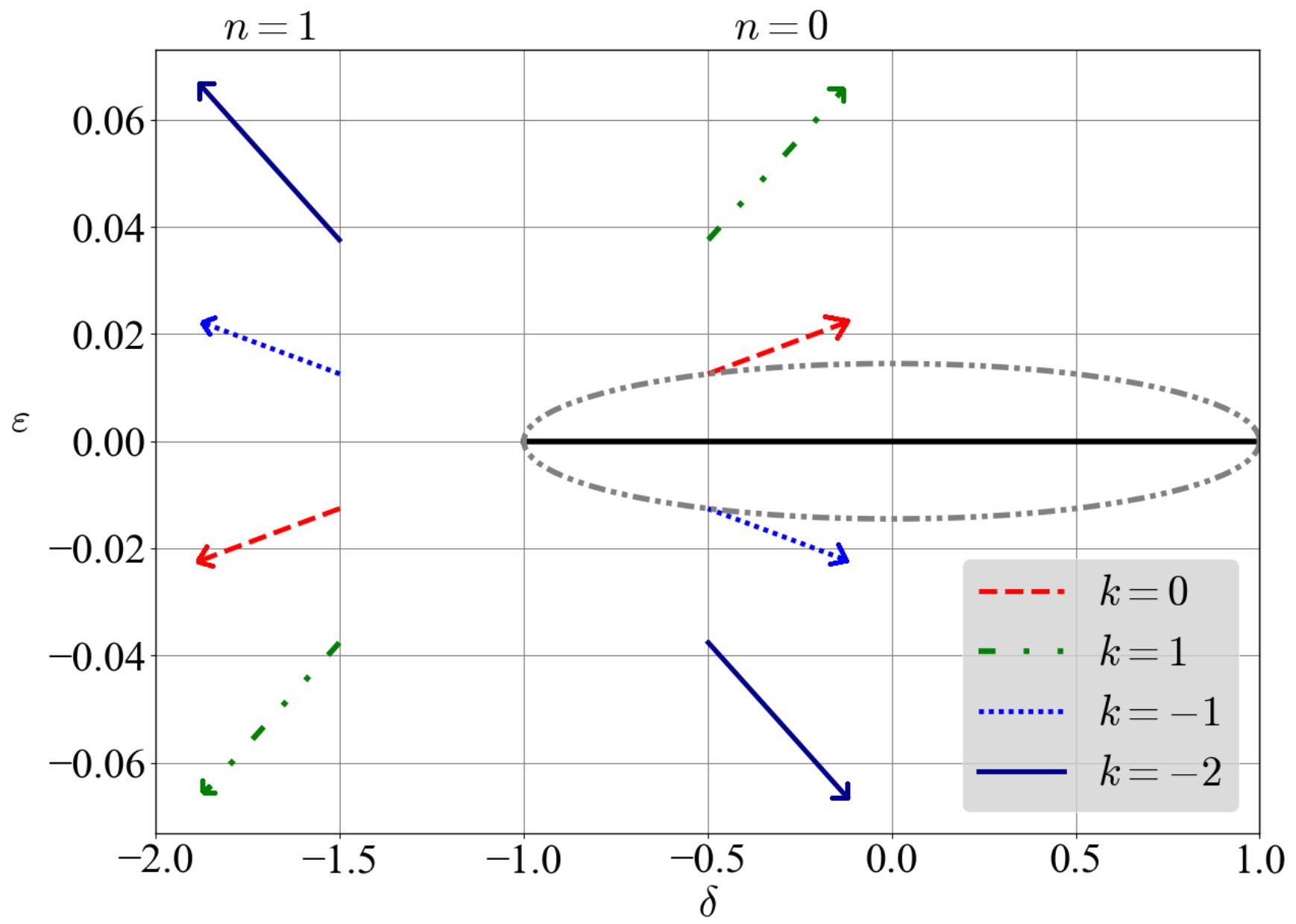} 
\caption{Segment $[-1,1]$ and the trajectories of singular points $z_{k,n}^{+}(t)\in\mathbb{C}$ of the exact solution to problem \eqref{eq_ub}, \eqref{eq_ub2} for $k = -2,-1,0,1,~n = 0,1$; $t \in[1,3.1]$ and $\nu = 0.001$. The Bernstein ellipse for the solution \eqref{exSol} at $t=1$ is plotted with a dash-dotted line.}
\label{sing_exSolDogan}
\end{figure}

Consider now the problem \eqref{eq_ub}, \eqref{eq_ub2}, which has the exact solution \eqref{exSol}. This solution has a countable number of simple poles \eqref{sing_exSol}. After transforming the square root of a complex number and making the change of variable $x$, which maps the segment $[0,2]$ onto $[-1,1]$, the coordinates of the singular points take the following form:
\begin{equation}
    \label{sing1_exSol}
    z_{k,n}^{\pm}=z_{k,n}^{\pm}(t)=-1\pm2|r_k(t)|^{1/2}\sqrt{\nu t} \exp\left(\dfrac{\vartheta_k(t)+2\pi n}{2}i \right)  ,~~k\in\mathbb{Z}~~,n\in \{0,1\},
\end{equation}
where $i=\sqrt{-1}$, $|r_k(t)| = \sqrt{\ln^2\sqrt{t_0/t} + \pi^2(2k+1)^2 }$, $\vartheta_k(t) = \arctan\dfrac{\pi(2k+1)}{\ln\sqrt{ t_0/t}}$ are the moduli and arguments of the radicands appearing in \eqref{sing_exSol}. The arguments are computed taking into account the quadrant, i.e. $\arctan(z)$ returns a value from the segment $[-\pi,\pi]$. We also note that $z_{k,0}^{+} = z_{k,1}^{-}$ and $z_{k,0}^{-} = z_{k,1}^{+}$, therefore, in what follows we can consider only $z_{k,n}^{+}$ for $n=0,1$. For $k= -1-j$ and $k = j$ with $j=0,1,2,...$ the poles \eqref{sing1_exSol} differ only by the sign of $\vartheta_k(t)$, which provides the symmetry of the poles with respect to the real axis. In Fig.~\ref{sing_exSolDogan} we show the domain of problem (the segment $[-1,1]$) and the trajectories of the singular points $z_{k,n}^{+}(t)$ in $\mathbb{C}$ for $k = -2,-1,0,1$, $n =0,1$, $t \in[1,3.1]$ and $\nu = 0.001$. From the figure, one can see that the singular points closest to the segment $[-1,1]$ are $z_{0,0}^{+}$ and $z_{-1,0}^{+}$. They are located symmetrically with respect to the real axis. In this case, Bernstein's theorem states that the error of interpolating the solution at $N$ Chebyshev nodes is proportional to $\rho^{-N}$, where $\rho$ is the sum of the semi-axes of an ellipse with foci at $\pm 1$ passing through the points $z_{0,0}^{+}$ and $z_{-1,0}^{+}$. This ellipse is called the Bernstein ellipse (see \cite{Tref}, chapter~8).

\begin{remark}
The presence of the poles $z_{0,0}^{+}$ and $z_{-1,0}^{+}$ has a significant impact on the gradients of the solution that can be observed in Fig.~\ref{fig1}b. These poles determine the rate of convergence of the proposed numerical method. In \cite{BE_CP}, for a certain class of problems for the BE having a countable number of poles, an influence of the whole set of singular points on the structure of the solution is studied, and the following conclusion is made: "the closest pole is typically one of primary interest in applications". Thus, computing the coordinates of the closest singular points with error control is the most important task.
\end{remark}

In Table~\ref{DoganvsUST}, the results of solving the problem \eqref{eq_ub}, \eqref{eq_ub2} by the proposed method are presented. For the cases $\nu = 0.005, 0.001$, a comparison with the results of Dogan et al. \cite{Dogan} was made. Here, No. denotes the number of the experiment. For computing the value $L_{\infty} = \|u_{ex}(T,x) - u(T,x)\|_\infty $ we use a grid of $10^4$ nodes uniformly distributed over the segment $[-1,1]$. The subscripts "$Dogan$" and "$STSM$" denote the results presented in \cite{Dogan} and the results of the STSM, respectively. The value of $err$ is defined by formula \eqref{Rerr}. The considerable advantage of the STSM is obvious. We note, however, that in \cite{Dogan} the parameters of the algorithm and the size of the grid, on which the value of $L_{\infty}$ $ _{Dogan} $ is computed, are not described. 

In the last columns of Table~\ref{DoganvsUST}, the errors of approximation for the coordinates of the singular point $z_{0,0}^{+}$ are presented: $$\delta_{err}=  \dfrac{\max\limits_{k=1,...,K}{|\delta_{ex}(\tau_k) - \delta(\tau_k) |}}{\max\limits_{k=1,...,K}{|\delta_{ex}(\tau_k)|}},~~\varepsilon_{err}= \dfrac{\max\limits_{k=1,...,K}{|\varepsilon_{ex}(\tau_k) - \varepsilon(\tau_k) |}}{\max\limits_{k=1,...,K}{|\varepsilon_{ex}(\tau_k)|}}.$$
Here, $\delta_{ex}$, $\varepsilon_{ex}$ and $\delta$, $\varepsilon$ are the exact and approximate coordinates of the singular point, respectively.

For the solutions with viscosities $\nu = 10^{-4}, 10^{-5}, 10^{-6}$, application of the composition $g\circ g$ of the form \eqref{gMap2} shows much better results. In Table~\ref{DoganvsUST} these solutions correspond to experiments No.~8, 10 and 11. For experiments No.~1--6, and also for experiment No.~9, the mapping \eqref{gMap} is used. In addition, when computing the solution for $\nu = 10^{-4}$ (experiment No.~7), to accelerate the algorithm we doubled the size $d_k$ of cells, whenever criterion \eqref{residd} with $n<200$ was satisfied while passing through 100 preceding cells (hereinafter, for brevity, we refer to this acceleration technique as "acceleration 100").

\begin{table}[h]
\caption{Results of the STSM algorithm for problem \eqref{eq_ub}, \eqref{eq_ub2} and comparison with the results of \cite{Dogan}}
\begin{center}
\begin{tabular}{l|l|l|l|l|l|l|l|l } 
         $\nu$ &     No.    & $T$ & $ L_{\infty}$ $ _{Dogan} $ & $ L_{\infty}$ $ _{STSM} $ & $err$ $ _{STSM} $&$N,M,K$   & $\delta_{err}$ & $\varepsilon_{err}$ \\
 \hline  
 \multirow{3}{*}{0.005} & 1 & 1.7 & 8.09e-03                   & 1.03e-11  & 3.27e-11 & 200, 10, 50    & 5.60e-07                        & 6.76e-06  \\
 \cline{2-9}            & 2 & 2.4 & 1.16e-02                   & 9.54e-12  & 3.64e-11 & 200, 10, 50    & 5.45e-07                        & 6.21e-06  \\
 \cline{2-9}            & 3 & 3.1 & 1.58e-02                   & 9.95e-12  & 4.35e-11 & 200, 10, 100   & 5.50e-07                        & 6.43e-06  \\ 
 \hline 
 \multirow{3}{*}{0.001} & 4 & 1.7 & 3.66e-02                   & 9.39e-11  & 2.57e-10 & 200, 10, 40    & 1.61e-05                        & 8.27e-05  \\
 \cline{2-9}            & 5 & 2.5 & 1.08e-02                   & 2.17e-11  & 7.25e-11 & 200, 10, 50    & 1.09e-07                        & 1.24e-05  \\
 \cline{2-9}            & 6 & 3.2 & 4.08e-02                   & 4.36e-11  & 1.64e-10 & 200, 10, 100   & 4.50e-07                        & 4.79e-05  \\
 \hline 
 \multirow{2}{*}{1.e-4} & 7 & 3.9 &                            & 7.32e-08  & 7.32e-08 & 200, 10, 488   & 2.25e-04                        & 1.59e-02  \\
 \cline{2-9}            & 8 & 3.9 &                            & 4.14e-09  & 1.65e-08 & 200, 10, 400   & 2.62e-04                        & 3.35e-02  \\
 \hline 
 \multirow{2}{*}{1.e-5} & 9 & 3.9 &                            & 6.07e-05  & 2.39e-04 & 200, 10, 11600 & 7.15e-06                        & 2.17e-05  \\
 \cline{2-9}            & 10 & 3.9 &                            & 9.39e-08  & 3.71e-07 & 200, 10, 4000  & 2.59e-05                        & 6.37e-02  \\
 \hline 1.e-6           & 11 & 3.9 &                            & 3.27e-06  & 1.29e-05 & 200, 10, 40000 & 1.21e-04                        & 6.87e-02         
 \end{tabular}
 \end{center}
\label{DoganvsUST}
\end{table} 

\begin{figure}[h!]
\centering 
\includegraphics[width=0.4\linewidth]{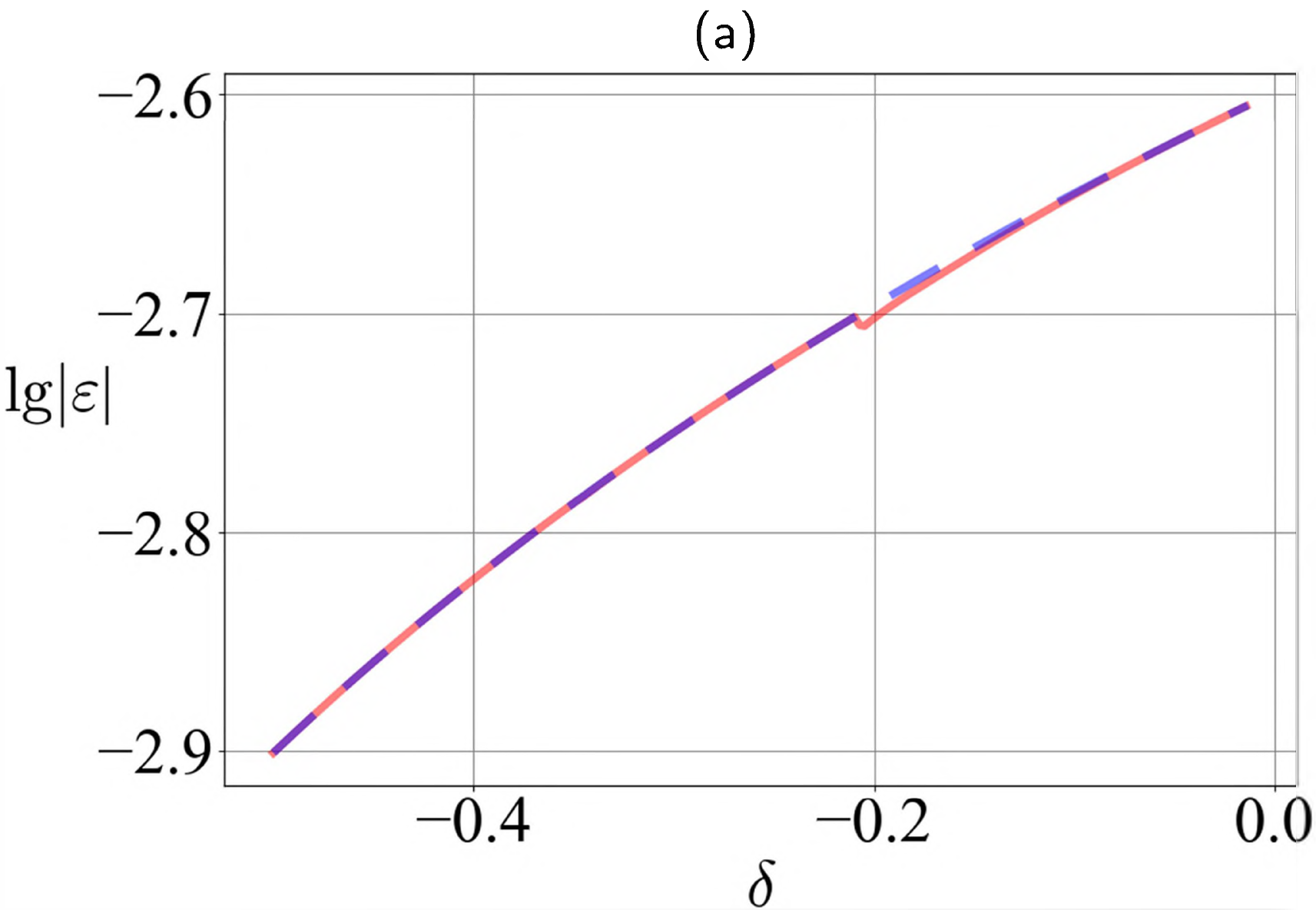}
\includegraphics[width=0.4\linewidth]{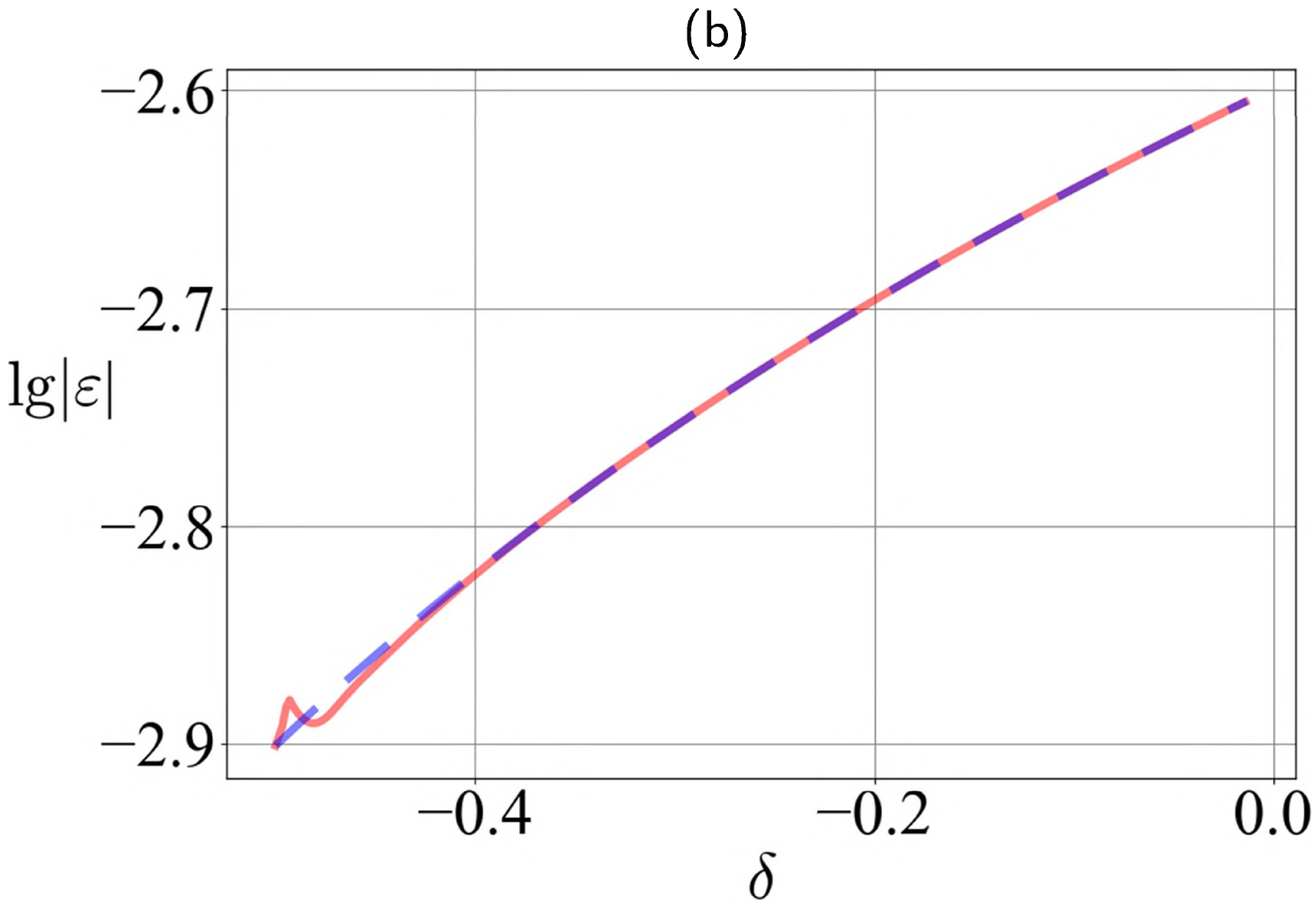}
\caption{Trajectories of the singular points of the exact solution (blue dashed line) and the numerical solution (red solid line) for (a) $\nu = 10^{-4}$ and $N = 200$, $M=10$, $K=488$ computed using the mapping \eqref{gMap} and the "acceleration 100" with initial size of a cell  $d_0 = 0.001$; and for (b) $\nu = 10^{-4}$ and $N = 200$, $M=10$, $K=400$ computed using the mapping \eqref{gMap2} without acceleration.}
\label{singDogan}
\end{figure}

In Fig.~\ref{singDogan}, we present the exact trajectory of the singular point $z_{0,0}^{+}$ (dashed line) and the approximate trajectory (solid line) obtained in computations (a) using the mapping \eqref{gMap} and the "acceleration 100", and (b) using the mapping \eqref{gMap2} without acceleration. The errors of computing the coordinates of the singular point for these regimes are presented in Table~\ref{DoganvsUST}, see experiments No.~7 and 8, respectively. It can be seen that the approximate trajectory agrees with the exact one with high accuracy; however, for both of the mentioned regimes (a) and (b), there are small deviations of the trajectories. Moreover, in case (a) such deviations can occur at the moment of time when the size of a cell is doubled, while in case (b) they can occur when the distance from the singular point to the problem domain becomes pretty small. It should be noted that increasing the size of cell can lead, on the one hand, to both an acceleration of the computations and an increase of the accuracy, and, on the other hand, to the appearance of jumps and oscillations of the approximate trajectory of the singular point. Similar results are obtained using mapping \eqref{gMap2} instead of \eqref{gMap}. However, in this case, decreasing $K$ is optional for accelerating the algorithm, whereas it is necessary for improving accuracy. Therefore, for the experiments No.~8, 10 and 11, the minimal value of $K$, for which the algorithm converges, was found experimentally. In this case, the error of the solution is much smaller than for larger $K$.

\begin{figure}[h!]
\centering  
\includegraphics[width=0.55\linewidth]{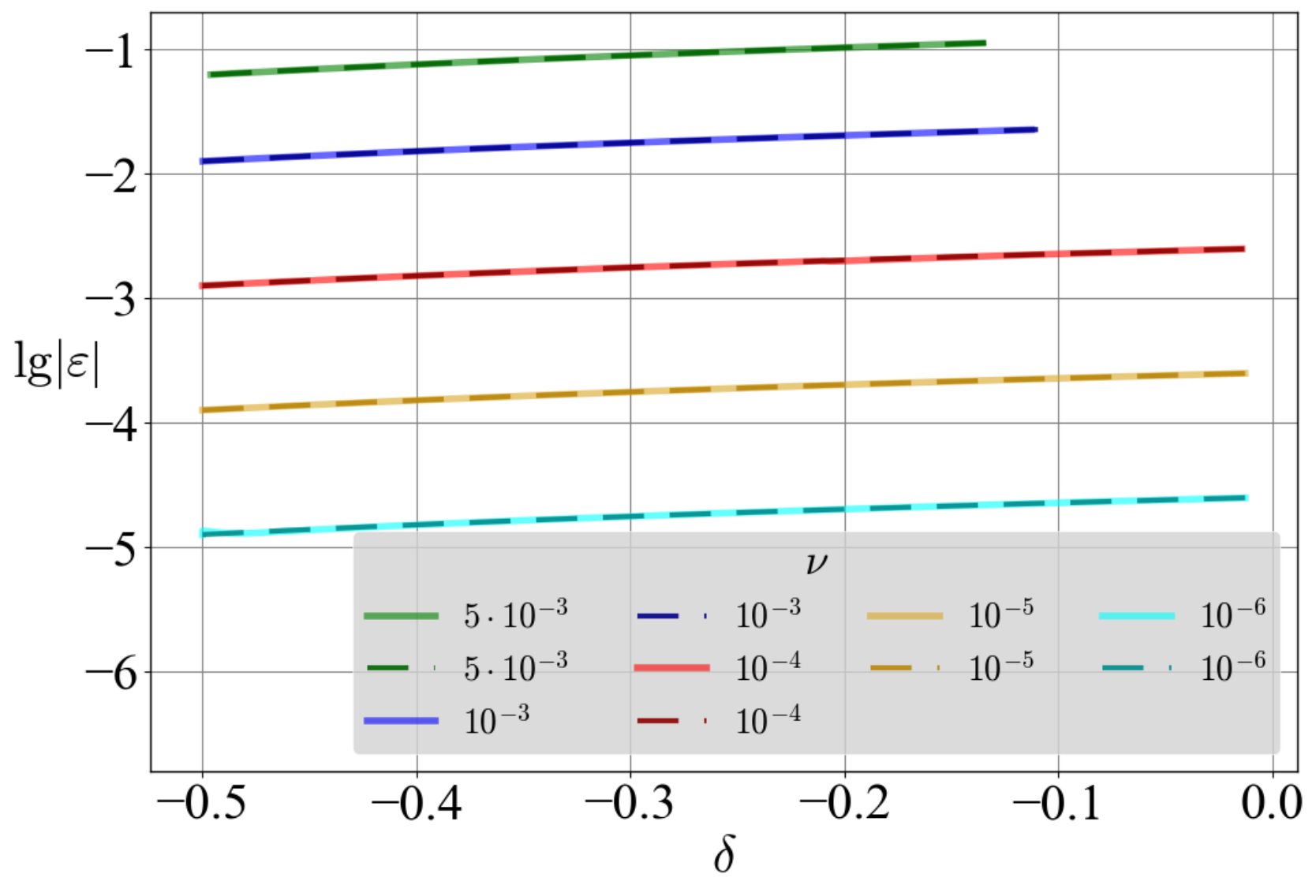}
\caption{ Approximate trajectories of the singular point $z_{0,0}^{+}(t)$ (solid lines) and exact trajectories (darker dashed lines) for $t\in[1,T]$ (the values of $T$ are given in Table~\ref{DoganvsUST}) and for various $\nu$.}
\label{singDogan1}
\end{figure}

In Fig.~\ref{singDogan1}, the approximate trajectories of the singular point $z_{0,0}^{+}(t)$ (solid lines) and the exact trajectories (darker dashed lines) are presented for the experimental setups from Table~\ref{DoganvsUST}. The trajectories are plotted for $t\in[1,T]$ (the values of $T$ are given in the table). From the figure, it can be seen that for all values of $\nu$ the constructed approximations of the trajectories have high accuracy.

\section{Concluding remarks}

\indent\indent We conclude our study by summarizing the most important properties of the proposed method. It enables one to solve numerically one dimensional problems of fluid dynamics with extremely small value of the viscosity. When solving initial boundary-value problems for the Burgers equation (BE), it becomes possible to perform the computations in floating-point arithmetic for the values of the viscosity $\nu$ ranging from $10^{-5}$ down to $10^{-8}$ and to systematically investigate the error. To the best of our knowledge, no such results have been reported in the literature yet.

In this context we should only refer to the paper \cite{Shap}, in which the solution to the BE with the initial data $u(0,x)=\varphi(x)=\cos(\pi x/2)$ and the boundary conditions $u(t,0)=1$, $u(t,2)=-1$ was computed for $x\in[0,2]$, $t>0$ by using a collocation and least squares method (CLSM) with $\nu$ reduced down to $\nu=10^{-14}$ and $\nu=0$. In this paper equation \eqref{BE} was written in a divergence form, and to approximate the solution with respect to the spatial variable, the CLSM was applied. It is known that the application of the least squares method minimizes the residual of a system of linear equations arising from the discretization of the differential equation. This suppresses the accumulation of round-off errors and yields a more stable algorithm. The obtained solutions do not have oscillations and converge to the solution of the Hopf equation as $\nu\to0$, which confirms the idea from \cite{Hopf}. However, the comparison with the exact solutions of the problem and the study of convergence of the numerical solutions for the fixed values of $\nu$ are absent in \cite{Shap}.

The key aspects of the method proposed in our paper, from the point of view of solving the problems with small $\nu$ are the application of rational approximations; the accounting for the location of the pole of the analytical continuation of the solution that is closest to the problem domain; the application of the new mappings \eqref{gcg}, \eqref{gMap2} for considerably expanding the domain of analyticity of the desired solution; and, as well, the proposed approach to forming and solving the linear algebra problem. In the last aspect the important observation is that the condition numbers of the matrices $R_A$ and $R_B$ presenting in the spectral decompositions \eqref{eq_spectral} are small, see \cite{SMNF, SemDiss} for the details. We should stress that despite the availability of formulas for the exact solution \eqref{eq_ub_exact}, verifying numerical solutions for small viscosity is highly non-trivial: the calculation of the exact solution with the error control requires a modification of the formulas for exact solution. We studied this problem in section~\ref{S2}.

It is worth noting that the discussed algorithms due to possibilities of fine-tuning their parameters (see subsection~\ref{secc2}) and of combining the rational barycentric interpolations with a high-order Runge--Kutta method, allow one to achieve a good balance between accuracy and computational costs, which is crucial for many applied problems.

The proposed method allows one to obtain a highly accurate approximation of the trajectories of singular points, even when they come extremely close to the domain of the problem. Hence, it is a promising tool for further investigation of the formation and evolution of singularities in hydrodynamic equations. It is well known that the emergence of singular points in the analytic continuations of solutions to the corresponding problems is related to a wide range of phenomena that are important for applications \cite{FrishBec, Calfisch, Weidemann1, Weidemann2, Zhuravleva}. Moreover, the motion of a singularity towards the problem domain is closely related to the question of the existence in time of smooth solutions of hydrodynamic equations which constitutes a famous problem. The design of new approaches to approximating and finding the corresponding solutions will open new ways for investigating this problem.

\section*{Acknowledgements}
The work of I.A. Bugoets is supported by the Mathematical Center in Akademgorodok under the agreement No.~075-15-2025-348 with the Ministry of Science and Higher Education of the Russian Federation. The authors are grateful to Vasilii Pavlovich Shapeev for providing a literature overview on numerical methods for the Burgers equation and for valuable discussions. The authors would like to thank Luka Sergeevich Bryndin for the provided software package and for discussing the implementation of the method.
 
\bibliography{cas-refs}
\end{document}